\documentclass[11pt]{article}
\usepackage{xypic}
\usepackage{amssymb}
\usepackage{amsfonts}
\usepackage{amsmath}
\usepackage[mathscr]{eucal}
\usepackage{latexsym}
\usepackage[french]{babel}

\newtheorem{teor}{Th{\'e}or{\`e}me}[section]          
\newtheorem{prop}[teor]{Proposition}
\newtheorem{cor}[teor]{Corollaire}                                            
 \newtheorem{rem}[teor]{Remarque}    
\newtheorem{lema}[teor]{Lemme}                                             
\newtheorem{defi}[teor]{D\'efinition}

\begin{document}
\author{Gentiana Danila\\
%%Institut de Math\'ematiques de Jussieu, 
%%UMR 7586 du CNRS, Case Postale 7012\\
%%2, Place Jussieu, 
%%75251 Paris Cedex 05, France\\
%%e-mail: gentiana@math.jussieu.fr\\
 Mathematical Institute, University of Warwick,\\
 Coventry CV4 7AL, United Kingdom\\
e-mail: gentiana@maths.warwick.ac.uk}
\title{Sur la cohomologie de la puissance sym\'etrique du fibr\'e tautologique sur le sch\'ema de Hilbert ponctuel d'une surface}
\date{1er avril 2001}
%\date{\ \ \ \ }
\maketitle

%%% espacements
\def\vs{\vskip 0.5cm}
%%%mode math

\def\has{{\rm H}}
\def\te{{\rm T}}
\def\be{{\rm B}}
\def\de{{\rm D}}
\def\tZ{\tilde{Z}}
\def\tmu{\tilde{\mu}}
\def\ttmu{\tilde{\tmu}}
\def\tpi{\tilde{\pi}}
\def\tnu{\tilde{\nu}}
\def\tp{\tilde{p}}

\def\ka{{\rm K}}
\def\es{{\rm S}}
\def\er{{\rm R}}
\def\supp{{\rm supp}\,}

\def\Hom{{\rm Hom}}
\def\Pic{{\rm Pic}\,}
\def\Tor{{\rm Tor}\,}
\def\Spec{{\rm Spec}\,}
\def\Ker{{\rm Ker}\,}
\def\Im{{\rm Im}\,}
\def\uTor{\underline{{\rm Tor}}\,}
\def\pgl{{\rm PGL}\,}
\def\gl{{\rm GL}\,}

%%%calligraphie

\def\maH{{\mathcal{H}}}
\def\maD{{\mathcal{D}}}
\def\maO{{\mathcal{O}}}
\def\maL{{\mathcal{L}}}
\def\maI{{\mathcal{I}}}
\def\sigm{{\mathfrak{S}}}
\def\sigmd{{\mathfrak{S}}_2}
\def\sigmn{{\mathfrak{S}}_n}
\def\emfrac{{\mathfrak{m}}}

%%%% espaces, groupes, ensembles

\def\proj{{\mathbb{P}}}
\def\pp{\proj_2}

\def\hilx{X^{\mbox{}^{[n]}}}
\def\tihilx{X_{\sim}^{\mbox{}^{[n]}}}
\def\trhilx{X_{(3)}^{\mbox{}^{[n]}}}
\def\ddhilx{X_{(2,2)}^{\mbox{}^{[n]}}}
\def\dhilx{\partial X^{\mbox{}^{[n]}}}
\def\hilxs{X_*^{\mbox{}^{[n]}}}
\def\hilxss{X_{**}^{\mbox{}^{[n]}}}
\def\dhilxs{\partial X_*^{\mbox{}^{[n]}}}
\def\hild{X^{\mbox{}^{[2]}}}
\def\ppd{\proj_2^{\mbox{}^{[2]}}}
\def\dhild{\partial X^{\mbox{}^{[2]}}}
\def\hilds{X_*^{\mbox{}^{[2]}}}
\def\dhilds{\partial X_*^{\mbox{}^{[2]}}}
\def\hilxnu{X^{\mbox{}^{[n,1]}}}
\def\hilxnd{X^{\mbox{}^{[n,2]}}}
\def\hilxndp{X^{\mbox{}^{[n,2]\prime}}}

\def\hilt{X^{\mbox{}^{[3]}}}
\def\hilts{X_*^{\mbox{}^{[3]}}}
\def\hiltss{X_{**}^{\mbox{}^{[3]}}}
\def\dhilt{\partial X^{\mbox{}^{[3]}}}
\def\dhilts{\partial X_*^{\mbox{}^{[3]}}}
\def\hiltc{X_{3C}^{\mbox{}^{[3]}}}
\def\hiltd{X^{\mbox{}^{[3,2]}}}
\def\xd{X^2}
\def\esxn{\es^nX}
\def\esxd{\es^2X}
\def\desxd{\partial\esxd}
\def\xn{X^n}
\def\xns{X^n_*}
\def\esdoi{\es^2}
\def\bedoi{\be^2}
\def\bedois{\be^2_*}
\def\benij{\be^n_{ij}}
\def\benijk{\be^n_{ijk}}
\def\benijkl{\be^n_{ij,kl}}
\def\ben{\be^n}
\def\den{\de^n}
\def\tiben{\be^n_{\sim}}
\def\tiden{\de^n_{\sim}}

\def\trben{\be^n_{(3)}}
\def\ddben{\be^n_{(2,2)}}
\def\tisigma{\Sigma_{\sim}}
\def\bens{\be^n_*}
\def\benss{\be^n_{**}}
\def\benud{\be^n_{12}}
\def\bends{\be^{\mbox{}^{[n,2]}}_*}
\def\bendij{\be^{\mbox{}^{[n,2]}}_{ij}}
\def\xij{\Xi_{ij}}
\def\xii{\Xi_{i}}
\def\xiud{\Xi_{[12]}}
\def\xiudi{\Xi_{[12]i}}
\def\xiund{\Xi_{12}}

\def\zed{{\mathbb Z}}
\def\comp{{\mathbb C}}

%%%fibr\'es, faisceaux

\def\dedoi{{\maD_2}}
\def\dna{{\maD_n^A}}
\def\dda{{\maD_2^A}}
\def\dta{{\maD_3^A}}
\def\dnl{{\maD_n^L}}
\def\ddl{{\maD_2^L}}
\def\dedoid{{\maD_2^{\mbox{}^{[2]}}}}
\def\dedoit{{\maD_2^{\mbox{}^{[3]}}}}
\def\dedoin{{\maD_2^{\mbox{}^{[n]}}}}

\def\elen{L^{\mbox{}^{[n]}}}
\def\elden{L^{{2\mbox{}^{[n]}}}}
\def\eled{L^{\mbox{}^{[2]}}}
\def\eltr{L^{\mbox{}^{[3]}}}
\def\eld{L^2}
\def\tield{\widetilde{L^2}}
\def\eldtr{L^{{2\mbox{}^{[3]}}}}
\def\eldoi{\maL_2}
\def\elded{L^{{2\mbox{}^{[2]}}}}
\def\elij{L_{ij}}

%%%fleches

\def\pnu{p_{n1}}
\def\pinu{\pi_{n1}}
\def\pnd{p_{n2}}
\def\pind{\pi_{n2}}
\def\pnus{p_{n1}^*}
\def\pinus{\pi_{n1*}}
\def\pnds{p_{n2}^*}
\def\pinds{\pi_{n2_*}}

\def\surto{\twoheadrightarrow}
\def\ra{\rightarrow}

%%%op\'erateurs

\def\tens{\otimes}

%%%abr\'eviations

\def\fs{{faisceau }}
\def\fx{{faisceaux }}
\def\alg{{alg\'ebrique }}
\def\algs{{alg\'ebriques }}
\def\th{{th\'eor\`eme }}

\def\thm{{th\'eor\`eme }}

\def\iso{{isomorphisme }}
\def\rep{{repr\'esentation }}
\def\reps{{repr\'esentations }}
\def\irr{{irr\'eductible }}
\def\irrs{{irr\'eductibles }}
\def\fib{{fibr\'e }}
\def\fibs{{fibr\'es }}
\def\mor{{morphisme }}
\def\mors{{morphismes }}
\def\sur{{surjectif }}
\def\co{{coh\'erent }}
\def\cohs{{coh\'erents }}
\def\app{{application }}
\def\apps{{applications }}

\def\resp{{respectivement }}
\def\inve{{inversible }}
\def\inves{{inversibles }}
\def\sct{{section }}
\def\scts{{sections }}
\def\cano{{canonique }}
\def\canos{{canoniques }}
\def\lin{{lin\'eaire }}
\def\lins{{lin\'eaires }}
\def\schil{{sch\'ema de Hilbert }}

{\bf Abstract~:} {\small We compute some cohomology spaces for the symmetric power of the 
tautological bundle tensor the determinant bundle on the punctual
Hilbert scheme $\hilx$ of a smooth projective surface $X$ on
$\comp$. }



Running heads: Fibr\'e tautologique sur le sch\'ema de Hilbert d'une surface

\section{Introduction}
\label{intro}

Soit $X$ une surface complexe projective et lisse et $L$ un fibr\'e inversible sur $X$. Pour tout entier $n$, on note $\hilx$ le sch\'ema de Hilbert qui param\`etre les sous-sch\'emas de $X$ de longueur $n$. Il est lisse et projectif de dimension $2n$ (\cite{Fogarty}). On suppose partout dans cet article que $n$ est un entier $\ge 2$.

On consid\`ere la vari\'et\'e d'incidence $\Xi=\hilxnu\subset\hilx\times X$ des points $(Z,x)$ qui v\'erifient $x\in \supp Z$.
On note $\pnu, \pinu$ les projections
\begin{equation}
\label{ecuatia00}
{\diagram
\Xi\rto^{\pnu}\dto_{\pinu}&X\\
\hilx.&
\enddiagram}
\end{equation}
On d\'efinit $\elen=\pinus(\pnus L)$. C'est un faisceau localement libre de rang $n$ sur $\hilx$. 
%
% On a r\'eduit dans \cite{Danila1} ce probl\`eme 
On s'int\'eresse au calcul des groupes de cohomologie:
$$\has^*(\hilx,\es^k\elen).$$
Le calcul est facile pour  $k=0$(\cite{Danila1}). Le calcul est compl\`etement r\'esolu dans \cite{Danila} pour $k=1$. Les r\'esultats sont r\'esum\'es par la formule:
\begin{equation}
\label{ecuatia002bis}
\has^*(\hilx,\es^k\elen)=\es^{n-k}\has^*(X,\maO)\otimes\es^k\has^*(X,L)\ \ {\rm\ \ pour\ \ }k=0,1,
\end{equation}
o\`u la puissance sym\'etrique est prise au sens $\zed/2$-gradu\'e.

\vs Le pr\'esent travail donne le calcul pour $k=2$ et  r\'esoud le probl\`eme pour toute la cohomologie pour  $n=2$ et $n=3$ et calcule l'espace vectoriel des sections de $\esdoi\elen $ pour  $n$ quelconque. Plus pr\'ecis\'ement, on introduit pour chaque entier $n$ l'application canonique
\begin{equation}
\label{ecuatia18}
can:\es^{n-2}\has^*(X,\maO_X)\otimes\esdoi\has^*(X,L)\to\has^*(\hilx,\esdoi\elen),
\end{equation}
%%%%%
%
%
%
%\begin{equation}
%\label{ecuatia2000}
%can:\es^{n-2}\has^*(X,A)\otimes\esdoi\has^*(X,L\otimes A)\to\has^*(\hilx,\esdoi\elen\otimes\dna)
%\end{equation}
de la mani\`ere suivante. On consid\`ere 
$$\esdoi_{\hilx}(\Xi)=(\Xi\times_{\hilx}\Xi)/{\sigmd}.$$ C'est un ferm\'e de $\hilx\times\esxd$. 
On note  $\pi, p $ les projections
\begin{equation}
\label{ecuatia44,5}
{\diagram
\esdoi_{\hilx}(\Xi)\rto^{p}\dto_{\pi}&\esxd\\
\hilx.&
\enddiagram}
\end{equation}
On a $\esdoi\elen=\pi_*p^*\ddl$. 
On note $P=(a,p)$ le \mor $\esdoi_{\hilx}(\Xi)\to \es^{n-2}X\times\esxd$, o\`u $a$ est le morphisme:
$$Z\mapsto (HC\circ\pi)(Z)-p(Z)\in \es^{n-2}X.$$
L'accouplement (\ref{ecuatia18}) r\'esulte de l'\iso \cano $\esdoi\elen=P^*(\maO\boxtimes\ddl))$, compte-tenu du fait que le \mor $\pi$ est fini.
%
%
%L'application canonique (\ref{ecuatia2000}) est la composition :
%%$$\displaylines{
%%\esdoi\has^*(X,L)=\has^*(\esxd,\dedoi)\stackrel{p^*}{\to}\hfill\cr
%%\hfill \has^*(\esdoi(\Xi),p^*\dedoi)=\has^*(\hilx,\pi_*p^*\dedoi)=\has^*(\hilx,\esdoi\elen).\cr}$$
%\begin{eqnarray*}
%\es^{n-2}\has^*(X,A)\otimes\esdoi\has^*(X,L\otimes A)=\has^*(\es^{n-2}X\times\esxd,\maD^A_{n-2}\boxtimes\maD^{L\otimes A}_2)\stackrel{P^*}{\to}\hfill\\
%\hfill\to\has^*(\esdoi_{\hilx}(\Xi),p^*\maD^L_2\otimes\pi^*\dna)=\has^*(\hilx,\pi_*p^*\maD^L_2\otimes\dna)=\has^*(\hilx,\esdoi\elen\otimes\dna).
%\end{eqnarray*}
%
Les r\'esultats principaux de cet article sont les suivants:
%
%
%Lorsque le faisceau \inve $A$ est trivial ($A=\maO$), le \mor canonique (\ref{ecuatia2000}) s'\'ecrit sous la forme:
%\begin{equation}
%\label{ecuatia18}
%can:\es^{n-2}\has^*(X,\maO_X)\otimes\esdoi\has^*(X,L)\to\has^*(\hilx,\esdoi\elen),
%\end{equation}
%et les th\'eor\`emes \ref{teor2.4'} et \ref{teor4.6'} peuvent s'\'ecrire plus  simplement:
\begin{teor}[n=2]
\label{teorema11}
%Soit $X$ une surface projective et lisse. Soit $L$ un fibr\'e inversible sur $X$. 
L'application canonique (\ref{ecuatia18})  pour $n=2$:
\begin{equation}
\label{ecuatia19}
can:\esdoi\has^*(X,L)\to\has^*(\hild,\esdoi\eled)
\end{equation}
induit la d\'ecomposition en somme directe:
$$\has^*(\hild,\esdoi\eled)=\esdoi\has^*(X,L)\bigoplus(\has^*(X,\maO_X)/\comp)\otimes\has^*(X,\eld).$$
\end{teor}
\begin{teor}[n=3]
\label{teor24}
%Soit $X$ une surface projective lisse et $L$ un faisceau inversible sur $X$. 
L'application canonique (\ref{ecuatia18})  pour $n=3$:
\begin{equation}
\label{ecuatia31}
can:\has^*(X,\maO_X)\otimes\esdoi\has^*(X,L)\to\has^*(\hilt,\esdoi\eltr),
\end{equation}
induit la d\'ecomposition en somme directe
$$\has^*(\hilt,\esdoi\eltr)=\has^*(X,\maO_X)\otimes\esdoi\has^*(X,L)\bigoplus(\esdoi \has^*(X,\maO_X)/{\has^*(X,\maO_X)})\otimes\has^*(X,\eld).$$
\end{teor}
Comme corollaire:
\begin{cor}
\label{coro}
Soit $X$ une surface projective lisse et $L$ un faisceau inversible sur $X$. Si $X$ satisfait en outre $q=p_g=0$, le morphisme canonique (\ref{ecuatia18}) est un isomorphisme $\has^*(\hilx,\esdoi\elen)=\esdoi\has^*(X,L)$ lorsque $n=2,3$.
\end{cor}

\begin{teor}[*=0]
\label{teorema154}
Le \mor canonique (\ref{ecuatia18}) est un \iso en degr\'e $0$:
$$
can:\esdoi\has^0(X,L)\stackrel{\sim}{\to}\has^0(\hilx,\esdoi\elen).
$$
\end{teor}
Les id\'ees des d\'emonstrations sont les suivantes: 
%Il n'y a aucune diff\'erence dans l'approche des cas $A$ \inve g\'en\'eral et $A=\maO_X$. Pour simplifier la pr\'esentation on se placera dans le cas $A=\maO_X$. 
%
Comprendre la structure du faisceau $\esdoi\elen$ \'equivaut \`a d\'ecrire le sch\'ema $\esdoi_{\hilx}(\Xi)$. Celui-ci a deux composantes irr\'eductibles:
\begin{eqnarray}
 &&\bullet \es_{(2)}=\{(Z,2x)\in\hilx\times\esxd,\ x\in\supp Z\}\nonumber\\
\label{ecuatia005}
&&\bullet X^{(n,2)}=\{(Z,x+y)\in\hilx\times\esxd,\ x+y\le HC(Z)\}.
\end{eqnarray}
Le probl\`eme se r\'eduit \`a l'\'etude de leur intersection sch\'ematique. Au-dessus de l'ouvert $\hilxss$ des sch\'emas lisses ces deux composantes sont disjointes. D\'ej\`a au-dessus de l'ouvert $\hilxs$ des sch\'emas avec au plus un point double cette intersection n'est pas triviale. Il est utile dans ce cas de faire le changement de base $\bens\to\hilxs$, o\`u $\ben$ est le produit fibr\'e:
\begin{equation}
\label{ecuatia900}
{\diagram
\ben\rto^{\rho}\dto_{q}&\xn\dto^{p}\\
\hilx\rto^{HC}&\esxn.
\enddiagram}
\end{equation}
et $\bens$ est l'image r\'eciproque $\bens=q^{-1}(\hilxs)$. 
Le probl\`eme se r\'eduit au cas $n=2$ (lemmes \ref{lema3}, \ref{lema5}). Un calcul de d\'eterminant (prop. \ref{prop2}) nous permet d'\'eviter le calcul de l'intersection sch\'ematique. 
Pour des raisons techniques on remplace $X^{(n,2)}$ avec le sch\'ema $\hilxnd\subset\hilx\times\hild$ param\'etrant les sous-sch\'emas $(Z,Z')$ tels que $Z\subset Z'$.
On note  $\pnd, \pind$ les projections
\begin{equation}
\label{ecuatia0}
{\diagram
\hilxnd\rto^{\pnd}\dto_{\pind}&\hild\\
\hilx&,
\enddiagram}
\end{equation}
et $\dedoin=\pinds(\pnds\dedoi)$.
On note $\dhilxs\subset\hilxs$ le ferm\'e des sch\'emas avec exactement un point double et $\mu:\dhilxs\to X$ l'application qui associe au sous-sch\'ema $Z\subset X$ son unique point double. Dans la section \ref{sectiune1} on prouvera:
\begin{teor}
\label{teorema1}
On a une suite exacte sur $\hilxs$:
\begin{equation}
\label{ecuatia000}
0\to\esdoi\elen\to\dedoin\oplus\elden\to \mu^*\eld\vert_{\dhilxs}\to 0.
\end{equation}
\end{teor}
Le \th \ref{teorema1} et le r\'esultat (\ref{ecuatia002bis}), (k=1), suffisent pour d\'emontrer \`a la section \ref{sectiune2} le \th \ref{teorema11}. Ce \th suffit pour d\'emontrer \`a la section \ref{sectiune5} le \th \ref{teorema154}.
Dans le cas $n > 2$ une \'etude fine des \mors $\pind:\hilxnd\to\hilx$ et $P:\hilxnd\to\es^{n-2}X\times \esxd$ est n\'ecessaire.
On r\'eussit \`a faire cette \'etude pour $n=3$, c'est l'objet de la section \ref{sectiune3}.
La section \ref{sectiune4} contient la preuve du \th \ref{teor24} \`a partir des r\'esultats de la section \ref{sectiune3}.
Dans la section \ref{sectiune6} on prolonge la suite exacte (\ref{ecuatia000})  de $\hilxs$ \`a une suite exacte sur un ouvert $\tihilx$ dont le compl\'ementaire est de codimension $3$ dans $\hilx$. En utilisant  ce fait on r\'esume dans la remarque \ref{prop67} ce qu'il reste \`a faire si on veut utiliser la m\^eme m\'ethode pour $n$ g\'en\'eral.

\vs Soit $A$ un \fib \inve sur $X$.
On consid\`ere le faisceau inversible $A\boxtimes \cdots\boxtimes A$ sur $\xn$. Le groupe $\sigmn$ agit sur $\xn$ par permutation des coordonn\'ees et cette action s'\'etend \`a une action \'equivariante sur $A\boxtimes\cdots\boxtimes A$. On d\'efinit le faisceau inversible 
\begin{equation}
\label{ecuatia001}
\dna=(A\boxtimes \cdots\boxtimes A)^{\sigmn}
\end{equation}
sur la vari\'et\'e $\esxn=\xn/{\sigmn}$. 
On note aussi $\dna$ l'image r\'eciproque de $\dna$ par le morphisme de Hilbert-Chow $HC:\hilx\to \esxn$, qui associe \`a un sous-sch\'ema $Z\subset X$ le cycle $\sum_{x\in X}{\rm lg}(Z_x)\cdot x$, o\`u $Z_x$ est la composante de $Z$ en $x$ et ${\rm lg}(Z_x)$ la longueur de $Z_x$.
La motivation de cet article a pour origine le calcul de l'espace de sections du fibr\'e d\'eterminant de Donaldson sur l'espace de modules de faisceaux semi-stables de rang $2$ sur le plan projectif (en analogie avec la  formule de Verlinde).
On a r\'eduit dans \cite{Danila1} ce probl\`eme au calcul des groupes de cohomologie:
$$\has^*(\hilx,\es^k\elen\otimes\dna).$$
L'analogue de la formule (\ref{ecuatia002bis}) pour $A$ inversible g\'en\'eral est (\cite{Danila1}, \cite{Danila}):
\begin{equation}
\label{ecuatia002}
\has^*(\hilx,\es^k\elen\otimes\dna)=\es^{n-k}\has^*(X,A)\otimes\es^k\has^*(X,L\otimes A)\ \ {\rm\ \ pour\ \ }k=0,1.
\end{equation}
L'analogue du \mor \cano (\ref{ecuatia18}) pour $A$ non trivial est:
\begin{equation}
\label{ecuatia2000}
can:\es^{n-2}\has^*(X,A)\otimes\esdoi\has^*(X,L\otimes A)\to\has^*(\hilx,\esdoi\elen\otimes\dna).
\end{equation}
Les r\'esultats qu'on prouve dans le cas $A$ g\'en\'eral sont:

\begin{teor}[n=2]
\label{teor2.4'}
%Soit $X$ une surface projective et lisse. Soient $L$ et $A$ deux fibr\'es inversibles sur $X$. 
L'application canonique (\ref{ecuatia2000}) pour $n=2$ fournit la d\'ecomposition canonique
$$\has^*(\hild,\esdoi\eled\otimes\dda)\simeq\esdoi\has^*(X,L\otimes A)\oplus K^*,$$
et $K^*$ rentre dans une suite exacte longue:
\begin{equation}
\label{ecuatia2001}
\cdots\to K^*\to\has^*(X,A)\otimes\has^*(X,\eld\otimes A)\to \has^*(X,\eld\otimes A^2)\to K^{*+1}\to\cdots
\end{equation}
\end{teor}

\begin{teor}[n=3]
\label{teor4.6'}
%Soit $X$ une surface projective et lisse. Soient $L$ et $A$ deux fibr\'es inversibles sur $X$. 
L'application canonique (\ref{ecuatia2000}) pour $n=3$ fournit la d\'ecomposition canonique
$$\has^*(\hilt,\esdoi\eltr\otimes\dta)\simeq\has^*(X,A)\otimes\esdoi\has^*(X,L\otimes A)\oplus K^*,$$
et $K^*$ rentre dans une suite exacte longue:
\begin{equation}
\label{ecuatia3001b}
\cdots\to K^*\to\esdoi\has^*(X,A)\otimes\has^*(X,\eld\otimes A)\to\has^*(X,A)\otimes\has^*(X,\eld\otimes A^2)\to K^{*+1}\to\cdots
\end{equation}

\end{teor}

\begin{teor}[*=0]
\label{teor154'}
Soit $n\ge 2$. L'application canonique (\ref{ecuatia2000})  fournit la d\'ecomposition canonique
\begin{equation}
\label{ecuatia2510}
\has^0(\hilx,\esdoi\elen\otimes\dna)=\es^{n-2}\has^0(X,A)\otimes\esdoi\has^0(X,L\otimes A)\oplus K_0,
\end{equation}
o\`u $K_0$ est le noyau du morphisme:
\begin{equation}
\label{ecuatia2515}
\es^{n-1}\has^0(X,A)\otimes\has^0(X,\eld\otimes A)\to\es^{n-2}\has^0(X,A)\otimes\has^0(X,\eld\otimes A^2),
\end{equation}
donn\'e par
\begin{equation}
\label{ecuatia2516}
u^{n-1}\otimes\alpha\mapsto(n-1)u^{n-2}\otimes u\alpha.
\end{equation}
\end{teor}
%Le \th \ref{teorema1} et le r\'esultat (\ref{ecuatia002}), (k=1), suffisent pour d\'emontrer \`a la section \ref{sectiune2} le \th \ref{teor2.4'}. Ce \th suffit pour d\'emontrer \`a la section \ref{sectiune5} le \th \ref{teor154}.
%
%Dans le cas $n\ge 2$ une \'etude fine des \mors $\pind:\hilxnd\to\hilx$ et $P:\hilxnd\to\es^{n-2}X\times \esxd$ est n\'ecessaire.
%On r\'eussit \`a faire cette \'etude pour $n=3$, c'est l'objet de la section \ref{sectiune3}.
%
%La section \ref{sectiune4} contient la preuve du \th \ref{teor4.6'} \`a partir des r\'esultats de la section \ref{sectiune3}.
%
%Dans la section \ref{sectiune6} on prolonge la suite exacte (\ref{ecuatia000})  de $\hilxs$ \`a une suite exacte sur un ouvert $\tihilx$ dont le compl\'ementaire est de codimension $3$ dans $\hilx$. En utilisant  ce fait on r\'esume dans la remarque \ref{prop67} ce qu'il reste \`a faire si on veut utiliser la m\^eme m\'ethode pour $n$ g\'en\'eral.
%

\vs {\bf Remerciements:} Je remercie M. Brion pour m'avoir attir\'e l'attention sur le lien entre les composantes (\ref{ecuatia005}) et les calculs de \cite{Danila1}. Je suis reconnaissante \`a J. Le Potier pour m'avoir propos\'e ce sujet et pour les id\'ees qu'il a g\'en\'ereusement partag\'e avec moi au cours de la r\'edaction de ce travail, en particulier le \th \ref{teorema13} et le lemme \ref{lema16}. Pendant sa r\'ealisation je me suis rejouie de l'ambiance d\'etendue de l'Institut de Math\'ematiques de l'Universit\'e de Warwick.

\section{ Suites exactes sur $\hilxs$}
\label{sectiune1}

On suppose partout dans cette section que $n$ est un entier $\ge 2$.
On fait la convention $\dedoi=\maD_2^L$ qui sera valable partout dans la suite. Le but de cette section est de d\'emontrer le \th \ref{teorema1}.
La preuve du \thm  utilise la proposition technique \ref{prop4} qui suit. On commence par d\'emontrer le cas $n=2$ dans la proposition suivante.
On note toujours $\eled$ le faisceau $q^*\eled$ sur $\bedoi$.
On note $E=q^{-1}(\dhild).$ C'est un diviseur sur $\bedoi$.

\begin{prop}
\label{prop2}
On a sur $\bedoi$ un diagramme commutatif:
$$\diagram
0\rto&\dedoi(-2E)\dto\rto&\esdoi\eled\dto\rto&\elded\dto\rto&0\\
0\rto&\dedoi(-E)\rto&\dedoi\rto&\dedoi\vert_{E}\rto&0,
\enddiagram$$
o\`u les suites horizontales sont exactes.
\end{prop}

%%%On note $\dhilx$ le ferm\'e de $\hilx$ constitu\'e par les sous-sch\'emas avec au moins un point double. Ce ferm\'e est un diviseur de $\hilx$.
%%%
%%%On note $\hilxs$ l'ouvert de $\hilx$ form\'e des sous-sch\'emas avec au plus un point double.
%%%
%%%Ol'intersection $\hilxs\cap\dhilx$.
%%%
%%%
%%%On consid\`ere le faisceau inversible $L\boxtimes L$ sur $\xd$. Le groupe $\sigmd$ agit sur $\xd$ par permutation des coordonn\'ees et cette action s'\'etend \`a une action \'equivariante sur $L\boxtimes L$. On d\'efinit le faisceau inversible $\dedoi=(L\boxtimes L)^{\sigmd}$ sur la vari\'et\'e $\esxd=\xd/_{\sigmd}$. On note aussi $\dedoi$ l'image r\'eciproque de $\dedoi$ par le morphisme de Hilbert-Chow $HC:\hild\to \esxd$.
%%%
%%%On consid\`ere la vari\'et\'e d'incidence $\hilxnd\subset\hilx\times\hild$ des points $(Z,Z')$ tels que $Z'$ est un sous-sch\'ema de $Z$.
%%%
%%%On note $\bedoi$ le produit fibr\'e:
%%%$$\diagram
%%%\bedoi\rto^{\rho}\dto_{q}&\xd\dto^{p}\\
%%%\hild\rto^{HC}&\esxd.
%%%\enddiagram$$
%%%
%%%On note \'egalement $\dedoi$ l'image r\'eciproque du faisceau inversible $\dedoi$  sur $\bedoi$. De la m\^eme mani\`ere, 

\vs\par{\bf Preuve de la proposition \ref{prop2}:}

On commence  par d\'emontrer l'exactitude de la premi\`ere ligne. On note $\Xi=\be^{\mbox{}^{[2,1]}}\subset\bedoi\times X$ le ferm\'e des couples $(\widetilde{Z},x)$ tels que $x\in\supp q(\widetilde{Z})$. On note \'egalement
$$\diagram
\Xi\rto^{p_{21}}\dto_{\pi_{21}}&X\\
\bedoi&.
\enddiagram$$
les morphismes vers $\bedoi, X$. Par le changement de base $\bedoi\stackrel{q}{\to}\hild$ il r\'esulte un isomorphisme sur $\bedoi$:
\begin{equation}
\label{ecuatia0,4}
\eled=\pi_{21*}(p_{21}^*L).
\end{equation}

On consid\`ere le sch\'ema $\esdoi_{\bedoi}(\Xi)=(\Xi\times_{\bedoi}\Xi)/\sigmd$ et 
%
%produit tensoriel $\Xi\times_{\bedoi}\Xi$. C'est un ferm\'e dans $\bedoi\times X\times X$. Le groupe $\sigmd$ agit en permutant les facteurs. 
%
%On note   C'est un ferm\'e dans $\bedoi\times\esdoi X.$ On note 
%
$pr_1, pr_2$ les projections:
$$\diagram
\esdoi_{\bedoi}(\Xi)\rto^{pr_2}\dto_{pr_1}&\esdoi X\\
\bedoi&.
\enddiagram$$
%On a un isomorphisme sur $\bedoi$:
%$$\eled\boxtimes\eled=\bar{pr_{1*}}(\bar{pr_{2*}}(L\boxtimes L)).$$
%On consid\`ere les invariants de cette \'egalit\'e pour l'action de $\sigmd$. 
On a:
$$\esdoi\eled=pr_{1*}(pr_2^*\dedoi).$$
On consid\`ere le morphisme diagonal $\Xi\stackrel{j}{\to}\esdoi_{\bedoi}(\Xi).$
Puisque $\pi_{21}:\Xi\to\bedoi$ est un morphisme fini, il est propre. Par cons\'equent le morphisme diagonal $j$ est une immersion. 

Le diagramme commutatif:
$$\diagram
&\Xi\rto^{p_{21}}\ddlto_{\pi_{21}}\dto^{j}&X\dto^{diag}\\
&\esdoi_{\bedoi}(\Xi)\rto^{pr_2}\dto^{pr_1}&\esdoi X\\
\bedoi\rto^{id}&\bedoi&
\enddiagram$$
d\'emontre que $j^*pr_2^*\dedoi=p_{21}^*\eld.$ Il r\'esulte une surjection sur $\bedoi$:
\begin{equation}
\label{ecuatia0,6}
\esdoi\eled=pr_{1*}(pr_2^*\dedoi)\stackrel{surj}{\to}\elded=\pi_{21*}(p_{21}^*\eld)\to 0.
\end{equation}

C'est une surjection entre un faisceau localement libre de rang $3$ et un faisceau localement libre de rang $2$. Son noyau est un faisceau inversible qui co\"{\i}ncide avec son d\'eterminant. Alors
$$\Ker surj=\det\esdoi\eled\otimes(\det\elded)^{-1}.$$
L'exactitude de la premi\`ere ligne r\'esulte du lemme:

\begin{lema}
\label{lema2,5}
Soit $X$ une surface projective lisse et $L$ un faisceau inversible. Pour le faisceau localement libre $\eled$ sur $\bedoi$ d\'efini par la relation (\ref{ecuatia0,4}) on a 
$$\det\eled=\dedoi(-E).$$ 
\end{lema}

Effectivement cela nous donne:
$$\Ker surj=(\det\eled)^3\otimes(\det\elded)^{-1}=\maD_2^3(-3E)\otimes(\maD_{2}^{\eld}(-E))^{-1}=\maD_2^3(-3E)\otimes(\maD_{2}^{-2}(E)=\dedoi(-2E).$$

\vs\par{\bf Preuve du lemme \ref{lema2,5}:}

On a d\'emontr\'e dans \cite{Danila}(2.10) l'existence d'une  suite exacte courte:
\begin{equation}
\label{ecuatia0,7}
0\to\eled\to p_1^*L\oplus p_2^*L\to\rho^*L\vert_{E}\to 0,
\end{equation}
o\`u $p_i, i=1,2$ sont les morphismes $\bedoi\stackrel{\rho}{\to}\xd\stackrel{pr_i}{\to} X$, et $\rho^*L$ est un faisceau inversible le long du diviseur $E$. Il en d\'ecoule l'\'egalit\'e
$$\det\eled=p_1^*L\otimes p_2^*L(-E)=\dedoi(-E). \ \ \ \ \ \ \ \ \ \hfill\Box$$

\vs Pour construire le carr\'e droit du diagramme de l'\'enonc\'e, on consid\`ere le diagramme commutatif
\begin{equation}
\label{ecuatia0,8}
{\diagram
\esdoi(\Xi)&\Xi\lto_{j}\\
\bedoi\uto^{a}&E\uto_{b}\lto_{i}
\enddiagram}
\end{equation}
o\`u $i$ est l'inclusion canonique, $a$ est le morphisme 
$$\widetilde{Z}\to(\widetilde{Z},HC\circ q(\widetilde{Z}))\in\esdoi(\Xi)\subset\bedoi\times\esxd,$$
et $b$ sa restriction \`a $E\subset\bedoi$.

D'apr\`es la construction de ces morphismes on obtient
\begin{eqnarray*}
a^*(pr_2^*\dedoi)&=&\dedoi,\\
i^*a^*(pr_2^*\dedoi)&=&\dedoi\vert_E.
\end{eqnarray*}
Par cons\'equent le diagramme
\begin{equation}
\label{ecuatia0,9}
{\diagram
\esdoi\eled\rto\dto&\elded\rto\dto&0\\
\dedoi\rto&\dedoi\vert_E\rto&0
\enddiagram}
\end{equation}
est commutatif.
Si on consid\`ere le noyau des suites horizontales du diagramme (\ref{ecuatia0,9}) on obtient le diagramme de l'\'enonc\'e. \hfill $\Box$

\vs Avant de donner la proposition \ref{prop4} analogue de la proposition \ref{prop2} pour le cas $n$ g\'en\'eral, on a besoin de plusieurs notations.

On rappelle que $\bens$ est l'image r\'eciproque $q^{-1}(\hilxs)$. 
C'est un rev\^etement non-ramifi\'e de degr\'e $n!$  de $\hilxs$. Le groupe $\sigmn$ agit sur $\bens$ et  $\hilxs=\bens/_{\sigmn}$. On a $\bedois=\bedoi$.

On note $\benss$ l'ouvert de $\bens$ des points $Z$ tels que $\rho(Z)$ est un $n$-uplet \`a termes distincts.

Pour tous $1\le i<j\le n$ on note $E_{ij}\subset\bens$ l'image r\'eciproque $\rho^{-1}(\Delta_{i,j})$. C'est un diviseur. Les diviseurs $E_{ij}$ sont deux \`a deux disjoints dans $\bens$.

On note $\benij$ l'ouvert $\benss\cup E_{ij}=\bens\setminus\bigcup_{\{k,l\}\ne\{i,j\}}E_{kl}.$

\begin{lema}
\label{lema3}
Il existe des applications 
$$(r_{ij},s_{ij}):\benij\to\bedoi\times X^{n-2}$$
qui identifient $\benij$ \`a un ouvert dans $\bedoi\times X^{n-2}$.
\end{lema}

\vs\par {\bf Preuve:}

D'apr\`es la d\'efinition de $\bens$, un point $\widetilde{Z}$ de $\bens$ consiste en un sous-sch\'ema $Z$ de longueur $n$ de $X$ et un $n$-uplet $(x_1,\cdots,x_n)\in\xn$ tels que $HC(Z)=x_1+\cdots+x_n.$
On d\'efinit 
$$s_{ij}(\widetilde{Z})=(x_1,\cdots,{\check{x_i}},\cdots,{\check{x_j}},\cdots,x_n)\in X^{n-2}.$$

D'apr\`es la d\'efinition de $\benij$, les termes $x_1,\cdots,{\check{x_i}},\cdots,{\check{x_j}},\cdots,x_n$ sont deux \`a deux distincts, et distincts de $x_i$ et $x_j$. Alors le sch\'ema $\maO_Z$ s'\'ecrit
 $$\maO_{Z'}\oplus\maO_{x_1}\oplus\cdots\oplus{\check{\maO_{x_i}}}\oplus\cdots\oplus{\check{\maO_{x_j}}}\oplus\cdots\maO_{x_n},$$
pour un sch\'ema $Z'$ de longueur $2$ \`a support $x_i+x_j$.

On d\'efinit $r_{ij}(\widetilde{Z})=(Z',(x_i,x_j))$, point de $\bedoi$.
 L'application r\'eciproque associe \`a
$$((Z',(x_i,x_j)),(x_1,\cdots,\check{x_i},\cdots,\check{x_j},\cdots,x_n))\in\bedoi\times X^{n-2}$$
le point
$$(\maO_{Z'}\oplus\maO_{x_1}\oplus\cdots\oplus\check{\maO_{x_i}}\oplus\cdots\oplus\check{\maO_{x_j}}\oplus\cdots\maO_{x_n},(x_1,\cdots,x_n))\in\benij.\ \ \ \ \ \ \ \ \ \hfill\Box$$

\vs On note $\xij\subset\benij\times\hild$ le graphe de l'application $q\circ r_{ij}:\benij\to\bedoi\to\hild$.

On note $\bends\subset\bens\times\hild$ le ferm\'e des points $(\widetilde{Z},Z')\in\bens\times\hild$ tels que $Z'\subset q(\widetilde{Z})$. On note encore:
\begin{equation}
\label{ecuatia800}
{\diagram
\bends\rto^{\pnd}\dto_{\pind}&\hild\\
\bens&
\enddiagram}
\end{equation}
les projections.

Par changement de base $\bens\to\hilxs$ on a un isomorphisme sur $\bens$:
\begin{equation} 
\label{ecuatia1}
\dedoin=\pinds(\pnds\dedoi).
\end{equation}

Dans la g\'en\'eralisation de la proposition \ref{prop2}, $\dedoin$ jouera le r\^ole de $\dedoi$.
On a $\xij\subset\bendij$, o\`u $\bendij=\pi_{n2}^{-1}(\benij)$. Autrement dit on a une surjection:
\begin{equation}
\label{ecuatia2}
\maO_{\bendij}\to\maO_{\xij}\to 0
\end{equation}
sur $\benij\times\hild$. On tensorise avec $\pnds\dedoi$ et on projette sur $\benij$. On obtient, en tenant compte de la relation (\ref{ecuatia1}) et de la d\'efinition de $\xij$, une surjection
\begin{equation}
\label{ecuatia3}
\dedoin\to r_{ij}^*\dedoi\to 0.
\end{equation}

On tensorise le dernier terme par $\maO_{E_{ij}}=\maO_{\benij}/_{\maI_{E_{ij}}}$. On obtient la surjection:
\begin{equation}
\label{ecuatia3,1}
\dedoin\to r_{ij}^*\dedoi\otimes\maO_{E_{ij}}\to 0
\end{equation}
sur $\benij$. Puisque le dernier terme a son support sur $E_{ij}$, l'application s'\'etend \`a tout $\bens$. On consid\`ere la somme de ces morphismes
\begin{equation}
\label{ecuatia3,2}
\dedoin\to \sum_{i<j}r_{ij}^*\dedoi\otimes\maO_{E_{ij}}\to 0.
\end{equation}
On note $\dedoin(-1)$ le noyau de cette application. Il jouera le r\^ole de $\dedoi(-E)$ dans l'analogue de la proposition \ref{prop2}.

De la m\^eme fa\c{c}on, si on tensorise le dernier terme de la surjection (\ref{ecuatia3}) par $\maO_{\benij}/_{\maI_{E_{ij}}^2}$, et on fait la somme de tous ces morphismes, on obtient une surjection:
\begin{equation}
\label{ecuatia3,3}
\dedoin\to \sum_{i<j}r_{ij}^*\dedoi\otimes \maO_{\benij}/_{\maI_{E_{ij}}^2} \to 0.
\end{equation}
On note par $\dedoin(-2)$ son noyau. Il jouera le r\^ole de $\dedoi(-2E)$ dans l'analogue de la proposition \ref{prop2}.

%Une derni\`ere notation. L'application $q:\bens\to\hilxs$ envoie $E_{ij}$ en $\dhilxs$. On note $q_{ij}:E_{ij}\to\dhilxs$ sa restriction \`a $E_{ij}$.
On est en mesure d'\'enoncer:

\begin{prop}
\label{prop4}
On a sur $\bens$ un diagramme commutatif:
\begin{equation}
\label{diag3,5}
{\diagram
0\rto&\dedoin(-2)\dto_{a}\rto&\esdoi\elen\dto\rto&\elden\dto\rto&0\\
0\rto&\dedoin(-1)\rto&\dedoin\rto&\sum_{i<j}r_{ij}^*\dedoi\vert_{E_{ij}}\rto&0,
\enddiagram}
\end{equation}
o\`u les suites horizontales sont exactes, et l'application $a$ est celle canonique. 

\end{prop}

\vs\par{\bf Preuve:}

L'id\'ee consiste \`a ramener le r\'esultat, \`a l'aide du lemme \ref{lema3}, au r\'esultat connu quand $n=2$. On d\'emontrera l'existence d'un diagramme avec les propri\'et\'es de l'\'enonc\'e pour chaque $\benij$ et on d\'emontrera qu'ils co\"{\i}ncident en restriction \`a $\benss$. Pour simplifier l'\'ecriture, on prendra $\benij=\benud$.

\begin{lema}
\label{lema5}
Pour $3\le i\le n$ et pour un faisceau inversible $L$ sur $X$, on note $L_i$ le faisceau $pr_i^*L$ sur $X^{n-2}$, o\`u $pr_i$ est la projection $X^{n-2}\to X$. Pour $3\le i\le n$ on note $\elij=L_i\otimes L_j$. Dans l'identification du lemme \ref{lema3} on a ($(i,j)=(1,2)$):

\begin{eqnarray}
\label{ecuatia4}
\elen&=&\eled\boxtimes\maO\oplus\maO\boxtimes(\sum_{3\le i\le n}L_i)\\
\label{ecuatia5}
\elden&=&\elded\boxtimes\maO\oplus\maO\boxtimes(\sum_{3\le i\le n}L_i^2)\\
\label{ecuatia6}
\esdoi\elen&=&\esdoi\eled\boxtimes\maO\oplus\eled\boxtimes(\sum_{3\le i\le n}L_i)\oplus\maO\boxtimes(\sum_{3\le i<j\le n }L_{ij})\oplus\maO\boxtimes(\sum_{3\le i\le n}L_i^2)\\
\label{ecuatia7}
\dedoin&=&\dedoi\boxtimes\maO\oplus\eled\boxtimes(\sum_{3\le i\le n}L_i)\oplus\maO\boxtimes(\sum_{3\le i<j\le n}L_{ij})\\
\label{ecuatia8}
\dedoin(-1)&=&\dedoi(-E)\boxtimes\maO\oplus\eled\boxtimes(\sum_{3\le i\le n}L_i)\oplus\maO\boxtimes(\sum_{3\le i<j\le n}L_{ij})\\
\label{ecuatia9}
\dedoin(-2)&=&\dedoi(-2E)\boxtimes\maO\oplus\eled\boxtimes(\sum_{3\le i\le n}L_i)\oplus\maO\boxtimes(\sum_{ 3\le i<j\le n}L_{ij})\\
\label{ecuatia10,1}
r^*_{12}\dedoi\vert_{E_{12}}&=&\dedoi\vert_{E}\boxtimes\maO\\
\label{ecuatia10,2}
r^*_{ij}\dedoi\vert_{E_{ij}}&=&0\ \ \ \ \mbox{\ \ \rm pour\ \ } \ \ \ (i,j)\ne (1,2).
\end{eqnarray}
\end{lema} 

En utilisant le lemme \ref{lema5}, le diagramme (\ref{diag3,5}) s'obtient sur $\benud$ de la somme directe des diagrammes :

$$\diagram
0\rto&\dedoi(-2E)\boxtimes\maO\dto\rto&\esdoi\eld\boxtimes\maO\dto\rto&\elded\boxtimes\maO\dto\rto&0\\
0\rto&\dedoi(-E)\boxtimes\maO\rto&\dedoi\boxtimes\maO\rto&\dedoi\vert_{E}\boxtimes\maO\rto&0,
\enddiagram$$

$$\diagram
0\rto&\eled\boxtimes(\sum_{i}L_i)\dto\rto&\eled\boxtimes (\sum_{i}L_i))\dto\rto&0\dto\rto&0\\
0\rto&\eled\boxtimes(\sum_{i}L_i)\rto&\eled\boxtimes (\sum_{i}L_i))\rto&0\rto&0,
\enddiagram$$
\begin{equation}
\label{diag10}
{\diagram
0\rto&\maO\boxtimes(\sum_{ i<j}L_{ij}) \dto\rto&\maO\boxtimes(\sum_{i<j}L_{ij})\oplus\maO\boxtimes(\sum_{i}L_i^2)\dto\rto&\maO\boxtimes(\sum_{i}L_i^2)\dto\rto&0\\
0\rto&\maO\boxtimes(\sum_{i<j}L_{ij})\rto&\maO\boxtimes(\sum_{i<j}L_{ij})\rto& 0\rto&0,
\enddiagram}
\end{equation}

o\`u le premier r\'esulte de la proposition \ref{prop2} et les deux autres sont \'evidents.

Le lemme suivant montre que la restriction du diagramme (\ref{diag3,5}) de $\benud$ \`a $\benss$ ne d\'epend pas du choix du couple $(1,2)$ parmi les couples $(i,j)$ avec $1\le i<j\le n$.

\begin{lema}
\label{lema6}
On consid\`ere pour $i=1,2$, les faisceaux inversibles $L_i=pr_i^*L$ sur $\xd$, o\`u $pr_i:\xd\to X$ est la projection sur la $i$-\`eme composante. On note $L_{12}=L_1\otimes L_2$. On a, dans l'identification $\bedoi\setminus E=\xd\setminus\Delta$, les isomorphismes:
\begin{flushleft}
\begin{eqnarray*}
\eled&=&L_1\oplus L_2,\\
\elded&=&L_1^2 \oplus L_2^2,\\
\esdoi\eled&=&L_1^2 \oplus L_{12}\oplus L_2^2,\\
\dedoi&=&\dedoi(-E)=\dedoi(-2E)=L_{12},\\
\dedoi\vert_{E}&=&0 \ \ \mbox{\ \ \rm sur \ \ }\ \ \bedoi\setminus E=\xd\setminus \Delta.
\end{eqnarray*}
\end{flushleft}
\end{lema}
Dans l'isomorphisme du lemme \ref{lema3}, $\benss$ s'identifie \`a un ouvert dans $(\bedoi\setminus E)\times X^{n-2}=(\xd\setminus\Delta)\times X^{n-2}$. Du lemme \ref{lema6} il r\'esulte que dans cette identification le diagramme (\ref{diag3,5}) se restreint sur $\benss$ au diagramme canonique (\ref{diag10}), o\`u les indices $i$ parcourent $1\le i\le n$ et les indices $i,j$ parcourent $1\le i<j\le n$. Par suite la restriction du diagramme (\ref{diag3,5}) de $\benij$ \`a $\benss$ ne d\'epend pas du choix du couple $(i,j)$. Par cons\'equent l'application $a$ de l'\'enonc\'e est l'application canonique. \hfill $\Box$

\vs\par {\bf Preuve du lemme \ref{lema5}:}

On se place dans les notations de la preuve du lemme \ref{lema3} pour $(i,j)=(1,2)$.
On rappelle que $\Xi\subset\bens\times X$ est le ferm\'e des points $(\widetilde{Z},x)$ pour lequels $x\in\supp q(\widetilde{Z})$. Puisque les points $(x_i)_{3\le i\le n}$ sont distincts deux \`a deux et disjoints du $\supp Z'$ le long de l'ouvert $\benud$, on obtient qu'au-dessus de cet ouvert le sch\'ema $\Xi$ s'\'ecrit comme la r\'eunion des sch\'emas disjoints $\xiud$ et $\xii$, $3\le i\le n$, o\`u:
\begin{eqnarray*}
\xiud&=&\{(\widetilde{Z},x), \widetilde{Z}\in\benud, x\in X, x\in\supp\tZ\},\\
\xii&=&\{(\tZ,x),\tZ\in\benud, x=x_i\}.
\end{eqnarray*}
Alors
\begin{equation}
\label{ecuatia9,1}
\elen=\pinus(\maO_{\Xi}\otimes \pnus L)=\pinus(\maO_{\xiud}\otimes \pnus L)\oplus\bigoplus_{3\le i\le n}\pinus(\maO_{\xii}\otimes \pnus L),
\end{equation}
o\`u $\pinu, \pnu$ sont les projections de $\benud\times\hild$ sur $\benud$ respectivement $\hild$.
Par d\'efinition $\xiud$ est le produit tensoriel de $\Xi\subset\bedoi\times X$ par l'application $r_{12}:\benud\to\bedoi.$
On obtient:
\begin{equation}
\label{ecuatia9,2}
\pinus(\maO_{\xiud}\otimes \pnus L)=r_{12}^*(\pinus(\maO_{\Xi}\otimes \pnus L))=r_{12}^*\eled.
\end{equation}
Par d\'efinition $\xii$ est le graphe de l'application $\benud\stackrel{\rho}{\to}\xn\stackrel{pr_i}{\to}X.$
On obtient:
\begin{equation}
\label{ecuatia9,3}
\pinus(\maO_{\xii}\otimes \pnus L)=L_i.
\end{equation}
Les relations (\ref{ecuatia9,1}), (\ref{ecuatia9,2}) et (\ref{ecuatia9,3}) impliquent la relation (\ref{ecuatia4}). Les relations (\ref{ecuatia5}) et (\ref{ecuatia6})
r\'esultent de la relation 
(\ref{ecuatia4}).

On rappelle que $\bends\subset\bens\times\hild$ est le ferm\'e des points $(\tZ,\xi)$ pour lequels le sch\'ema $\xi$ est inclus dans le sch\'ema $q(\tZ)$. Puisque les points $(x_i)_{3\le i\le n}$ sont distincts deux \`a deux et disjoints de $\supp Z'$ le long de $\benud$, on trouve qu'au-dessus de $\benud$ le sch\'ema $\bends$ est la r\'eunion des sch\'emas disjoints:
\begin{eqnarray*}
\xiund&=&\{(\tZ,Z'), \tZ\in\benud\},\\
\xiudi&=&\{(\tZ,\xi),  \tZ\in\benud, \xi=\maO_{x_1}\oplus\maO_{x_i}  \ \ \mbox{\rm o\`u\ \ } \xi= \maO_{x_2}\oplus\maO_{x_i}, 3\le i\le n \},\\
\xij&=&\{(\tZ,\xi),  \tZ\in\benud, \xi=\maO_{x_i}\oplus\maO_{x_j} \}.
\end{eqnarray*}
On rappelle qu'on a not\'e au (\ref{ecuatia800}) $\pind, \pnd$
 les projections de $\bends$ sur $\bens$ et respectivement $\hild$. La relation (\ref{ecuatia1}) implique
\begin{equation}
\label{ecuatia9,4}
\dedoin=\pinds(\maO_{\xiund}\otimes \pnds \dedoi)\oplus\bigoplus_{3\le i\le n}\pinds(\maO_{\xiudi}\otimes \pnds \dedoi)\oplus\bigoplus_{3\le i\le n}\pinds(\maO_{\xij}\otimes \pnds \dedoi).
\end{equation}
Par d\'efinition, $\xiund$ est le graphe de l'application $\benud\stackrel{r_{ij}}{\to}\bedoi\stackrel{\rho}{\to}\xd.$ On obtient:
\begin{equation}
\label{ecuatia9,5}
\pinds(\maO_{\xiund}\otimes \pnds \dedoi)=r_{12}^*\dedoi.
\end{equation}
Le m\^eme argument implique, pour $3\le i\le n$:
\begin{equation}
\label{ecuatia9,6}
\pinds(\maO_{\xij}\otimes \pnds \dedoi)=L_{ij}.
\end{equation}
D'apr\`es la d\'efinition de $\xiudi$, dans l'identification $\benud\subset\bedoi\times X^{n-2}$, le morphisme $\xiudi\to\benud$ est le produit des morphismes $\xiud\to\bedoi$ et $\xii\to X^{n-2}$:
$$\xiud\times\xii\subset(\bedoi\times X)\times(X^{n-2}\times X)=(\bedoi\times X^{n-2})\times\xd\to\bedoi\times X^{n-2}.$$
En outre, le morphisme $\xiudi\to\benud\times\hild\stackrel{pr_2}{\to}\hild$ co\"{\i}ncide avec le morphisme:
$$\xiud\times\xii\subset(\bedoi\times X^{n-2})\times(\xd\setminus\Delta)\to\xd\setminus\Delta\stackrel{p}{\to}\esdoi X\setminus\Delta\stackrel{\sim}{\leftarrow}\hild\setminus\dhild\subset\hild.$$
Puisque $p^*\dedoi=L\boxtimes L$ on trouve:
\begin{equation}
\label{ecuatia9,7}
\pinds(\maO_{\xiudi}\otimes \pnds \dedoi)=\pinds(\maO_{\xiud}\otimes \pnds L)\otimes \pinds(\maO_{\xii}\otimes \pnds L)=\eled\boxtimes L_i.
\end{equation}
Les relations (\ref{ecuatia9,4}), (\ref{ecuatia9,5}), (\ref{ecuatia9,6}) et (\ref{ecuatia9,7}) impliquent la relation 
(\ref{ecuatia7}).
Dans la notation de la relation (\ref{ecuatia7}), l'application (\ref{ecuatia3}) est la projection
\begin{equation}
\label{ecuatia9,75}
\dedoin\to\dedoi\boxtimes\maO.
\end{equation}
Dans l'identification $\benud\subset\bedoi\times X^{n-2}$, le diviseur $E_{12}$ co\"{\i}ncide avec $E\times X^{n-2}$. Par cons\'equent dans cette identification on a:
\begin{eqnarray*}
\maO_{E_{12}}&=&\maO_E\boxtimes\maO,\\
\maO_{\benud}/\maI^2_{E_{12}}&=&\maO_{\bedoi}/\maI^2_E\boxtimes\maO.
\end{eqnarray*}
On tensorise les deux derni\`eres relations par $r_{12}^*\dedoi$ et on obtient:
\begin{eqnarray}
\label{ecuatia9,8}
r_{12}^*\dedoi\vert_{E_{12}}&=&\dedoi\vert_{E}\boxtimes\maO,\\
\label{ecuatia9,85}
r_{12}^*\dedoi\otimes\maO_{\benud}/\maI^2_{E_{12}}&=&\dedoi\otimes\maO_{\bedoi}/\maI^2_E\boxtimes\maO.
\end{eqnarray}
La relation (\ref{ecuatia9,8}) est la relation (\ref{ecuatia10,1}) \'enonc\'ee.

Pour $(i,j)\ne(1,2)$ le faisceau $r_{ij}^*\dedoi\vert_{E_{ij}}$ est nul en restriction \`a $\benud$, son support \'etant sur $E_{ij}$. La relation (\ref{ecuatia10,2}) en d\'ecoule.

\`A partir des relations (\ref{ecuatia9,75}), (\ref{ecuatia9,8}) et (\ref{ecuatia10,2}), l'application 
(\ref{ecuatia3,2}) s'\'ecrit
$$\dedoi\boxtimes\maO\oplus\eled\boxtimes(\sum_{i}L_i)\oplus\maO\boxtimes(\sum_{3\le i\le n}L_{ij})\to\dedoi\vert_E\boxtimes\maO.$$
Par cons\'equent son noyau $\dedoi(-1)$ se met sous la forme  (\ref{ecuatia8}). De mani\`ere analogue, les relations (\ref{ecuatia9,75}), (\ref{ecuatia9,85}) et l'annulation du terme $r_{ij}^*\dedoi\otimes\maO_{\benij}/\maI^2_{E_{ij}}$ sur $\benud$ pour $\{i,j\}\ne\{1,2\}$ impliquent la relation (\ref{ecuatia9}).\hfill $\Box$

\vs\par {\bf Preuve du lemme \ref{lema6}:}

La suite exacte (\ref{ecuatia0,7}) implique $\eled=L_1\oplus L_2$ sur $\bedoi\setminus E$. Les deux identit\'es suivantes de l'\'enonc\'e en r\'esultent.
Par construction on  a $p^*\dedoi=L\boxtimes L=L_1\otimes L_2$ sur $\xd$. On obtient l'\'egalit\'e sur $\bedoi$:
$$\dedoi\stackrel{notation}{=}\rho^*p^*\dedoi=\rho^*(L_1\otimes L_2)=L_1\otimes L_2\stackrel{notation}{=}L_{12}.$$
Les \'egalit\'es 
$$\dedoi(-E)=\dedoi(-2E)=\dedoi$$
sur $\bedoi\setminus E$ et
$$\dedoi\vert_E=0$$
sur 
$\bedoi\setminus E$ sont \'evidentes.
\hfill $\Box$

\vs Cela termine la preuve de la proposition \ref{prop4} concernant l'existence du diagramme  (\ref{diag3,5}) sur $\bens$. Pour passer de la vari\'et\'e $\bens$  \`a  $\hilxs$ on utilise:

\begin{lema}
\label{lema7}

On a des isomorphismes de faisceaux sur $\hilxs=\bens/_{\sigmn}$:
\begin{flushleft}
\begin{eqnarray}
\label{ecuatia11}
(\dedoin(-1))^{\sigmn}&=&(\dedoin(-2))^{\sigmn}\\
\label{ecuatia12}
(\dedoin)^{\sigmn}&=&\dedoin\\
\label{ecuatia13}
(\esdoi\elen)^{\sigmn}&=&\esdoi\elen\\
\label{ecuatia14}
(\elden)^{\sigmn}&=&\elden\\
\label{ecuatia15}
(\sum_{i<j}r_{ij}^*\dedoi\vert_{E_{ij}})^{\sigmn}&=&\mu^*\eld\vert_{\dhilxs}.
\end{eqnarray}

\end{flushleft}

\end{lema}

En appliquant le lemme \ref{lema7} on obtient du diagramme (\ref{diag3,5}) le diagramme commutatif sur $\hilxs$:
\begin{equation}
\label{diag16}
{\diagram
0\rto&\dedoin(-2)^{\sigmn}\dto_{a}\rto&\esdoi\elen\dto\rto&\elden\dto\rto&0\\
0\rto&\dedoin(-1)^{\sigmn}\rto&\dedoin\rto&\mu^*\eld\vert_{\dhilxs}\rto&0,
\enddiagram}
\end{equation}
o\`u les lignes sont exactes et le morphisme $a$ est l'isomorphisme canonique.

\vs\par{\bf Preuve du \thm \ref{teorema1}:}

La suite exacte (\ref{ecuatia000}) est le mapping c\^one pour le diagramme (\ref{diag16}) en 
supprimant la colonne de gauche.
%tenant compte du fait que le morphisme $a$ est un isomorphisme. 
\hfill $\Box$

\vs\par{\bf Preuve du lemme \ref{lema7}:}

Les relations (\ref{ecuatia12}), (\ref{ecuatia13}), (\ref{ecuatia14}) sont \'evidentes, les faisceaux localement libres $\dedoin, \esdoi\elen$ et $\elden$ sur $\bens$ provenant de $\hilxs$.
On d\'emontrera la relation:
\begin{equation}
\label{ecuatia16,1}
\left(\sum_{i<j}r_{ij}^*\dedoi\otimes\maO_{\bens/\maI^2_{E_{ij}}}\right)^{\sigmn}=\mu^*\eld\vert_{\dhilxs}.
\end{equation}
La suite (\ref{ecuatia3,2}) qui d\'efinit $\dedoin(-1)$ est $\sigmn$-invariante. Si on consid\`ere la suite des $\sigmn$-invariants de cette suite et les relations (\ref{ecuatia12}) et (\ref{ecuatia15}) on trouve
\begin{equation}
\label{ecuatia16,2}
(\dedoin(-1))^{\sigmn}=\Ker(\dedoin\to \mu^*\eld\vert_{\dhilxs}).
\end{equation}
De la m\^eme mani\`ere, \`a partir des relations (\ref{ecuatia3,3}), (\ref{ecuatia12}) et (\ref{ecuatia16,1}) on obtient
\begin{equation}
\label{ecuatia16,3}
(\dedoin(-2))^{\sigmn}=\Ker(\dedoin\to \mu^*\eld\vert_{\dhilxs}).
\end{equation}
Les \'egalit\'es (\ref{ecuatia16,2}) et (\ref{ecuatia16,3}) impliquent l'\'egalit\'e (\ref{ecuatia11}).
Il reste \`a d\'emontrer (\ref{ecuatia15}) et (\ref{ecuatia16,1}). 

Le terme gauche de l'\'equation (\ref{ecuatia15}) est une somme directe de faisceaux param\'etr\'es par des sous-ensembles \`a deux \'el\'ements $\{i,j\}\subset\{1,2,\cdots,n\}$. Le groupe $\sigmn$ agit transitivement sur l'ensemble de ces indices et le stabilisateur du point $\{1,2\}$ est le sous-groupe $\sigmd\times\sigm_{n-2}$. D'apr\`es le lemme 2.2 de \cite{Danila} on obtient:
\begin{eqnarray}
(\sum_{i<j}r_{ij}^*\dedoi\vert_{E_{ij}})^{\sigmn}&=&(r_{12}^*\dedoi\vert_{E_{12}})^{\sigmd\times\sigm_{n-2}}\stackrel{a}{=}(\dedoi\vert_E\boxtimes\maO)^{\sigmd\times\sigm_{n-2}}=\nonumber\\
\label{ecuatia1000}
&=&(\dedoi\vert_E)^{\sigmd}\boxtimes\maO_{X^{n-2}}^{\sigm_{n-2}}
\end{eqnarray}
o\`u $a$ est l'identification du lemme \ref{lema3} pour $\{i,j\}=\{1,2\}$.
Un point de $\dhilxs$ s'\'ecrit $\maO_Z=\maO_{Z'}\oplus\maO_{x_3}\oplus\cdots\oplus\maO_{x_n}$ pour $Z'$ sch\'ema singulier de longueur $2$ et $x_i$ des points deux \`a deux disjoints et disjoints de $\supp Z'$. Alors $\dhilxs$ s'identifie avec un ouvert de $\dhild\times\es^{n-2}X$. Dans cette identification le \mor $\mu$ est la composition:
$$\dhilxs\to\dhild\times\es^{n-2}X\stackrel{pr_1}{\to}\dhild\stackrel{HC}{\to}X=\partial\esxd,$$
o\`u $HC$ est le \mor de Hilbert-Chow.
Comme $\maO_{X^{n-2}}^{\sigm_{n-2}}=\maO_{\es^{n-2}X}$, l'\'equation (\ref{ecuatia15}) se r\'eduit \`a l'\'equation suivante sur $\dhild$:
$$(\dedoi\vert_E)^{\sigmd}=HC^*\eld.$$
Par d\'efinition, le faisceau $\dedoi$ sur $\bedoi$ est l'image r\'eciproque $HC^*\dedoi$, o\`u $\dedoi$ est un faisceau inversible sur $\esxd$. On a $\dedoi\vert_{X=\partial\esxd}\simeq\eld$. Par cons\'equent l'\'egalit\'e (\ref{ecuatia15}) est \'equivalente \`a l'\'egalit\'e:
\begin{equation}
\label{ecuatia1001}
\maO^{\sigmd}_E=\maO_{\dhild}.
\end{equation}
De la m\^eme mani\`ere que l'\'equation (\ref{ecuatia1000}) on obtient:
$$\left(\sum_{i<j}r_{ij}^*\dedoi\otimes\maO_{\bens/\maI^2_{E_{ij}}}\right)^{\sigmn}=\left(\dedoi\otimes\maO_{\bedoi/\maI^2_{E}}\right)^{\sigmd}\boxtimes\maO_{X^{n-2}}^{\sigm_{n-2}}$$
et la relation (\ref{ecuatia16,1}) est \'equivalente \`a l'\'egalit\'e:
\begin{equation}
\label{ecuatia1002}
(\maO_{\bedoi/\maI^2_{E}})^{\sigmd}=\maO_{\dhild}.
\end{equation}
Les \'egalit\'es (\ref{ecuatia1001}), (\ref{ecuatia1002}) r\'esultent du fait que le \mor de degr\'e $2$, $q:\bedoi\to\hild=\bedoi/\sigmd$ est ramifi\'e au-dessus du diviseur $\dhild$: $q^*\maO_{\dhild}=\maO_{\bedoi/\maI^2_{E}}$.
\hfill $\Box$

\section{ Le calcul de $\has^*(\hild,\esdoi\eled)$}
\label{sectiune2}

On calculera dans cette section la cohomologie $\has^*(\hild,\esdoi\eled)$ \`a partir du \thm \ref{teorema1}. Dans le cas $n=2$ on a $\hilds=\hild$, $\dhilds=\dhild$  et $\dedoid=\dedoi$. Par cons\'equent le \thm \ref{teorema1} affirme l'existence d'une suite exacte:
\begin{equation}
\label{ecuatia17}
0\to\esdoi\eled\to\dedoi\oplus \elded\to\mu^*\eled\vert_{\dhild}\to 0
\end{equation}
sur $\hild$.

On commence par calculer la cohomologie sur $\hild$ des termes qui appara\^{\i}ssent dans la suite exacte ci-dessus.

\begin{prop}
\label{prop8}
On a
$$\has^*(\hild,\elded)=\has^*(X,\maO)\otimes\has^*(X,\eld).$$

\end{prop}

\vs\par{\bf Preuve:}

Ce r\'esultat est un corollaire de l'affirmation g\'en\'erale (\ref{ecuatia002bis}) pour $k=1$, $n=2$.
%%\begin{equation}
%%\label{ecuatia17,3}
%%\has^*(\hilx,\elen)=\es^{n-1}\has^*(X,\maO_X)\otimes\has^*(X,L),
%%\end{equation}
%%la puissance sym\'etrique \'etant prise au sens $\zed/2$-gradu\'e. 
\hfill $\Box$

\begin{lema}
\label{lema9}

On a 
$$\has^*(\hild,\dedoi)=\esdoi\has^*(X,L),$$ 
la puissance sym\'etrique \'etant prise au sens $\zed/2$-gradu\'e.
\end{lema}

\vs\par{\bf Preuve:}

Ce r\'esultat est un corollaire de l'affirmation g\'en\'erale (\ref{ecuatia002bis}) pour $k=0$, $n=2$.\hfill $\Box$

\begin{lema}\label{lema10}

On a
$$\has^*(\hild,\mu^*\eld\vert_{\dhild})=\has^*(X,\eld).$$
\end{lema}

\vs\par{\bf Preuve:}

On a $\has^*(\hild,\mu^*\eld\vert_{\dhild})=\has^*(\dhild,\mu^*\eld).$
Le morphisme $\mu$ est une fibration \`a fibres $\proj_1$. Alors
$$\begin{array}{ccc}
   \er^q\mu_*\maO_{\dhild}&=&
\begin{cases}
  0 \mbox{ \ \ \ si \ \ \ } q>0\cr
\maO_X \mbox{ \ \ \ si \ \ \ } q=0.
\end{cases}
 \end{array}$$
On en d\'eduit que
 $$\begin{array}{ccc}
   \er^q\mu_*(\mu^*\eld)&=&
\begin{cases}
  0 \mbox{ \ \ \ si \ \ \ } q>0\cr
\eld \mbox{ \ \ \ si \ \ \ } q=0.
\end{cases}
 \end{array}$$
Par la suite spectrale de Leray on obtient:
$$\has^*(\hild,\mu^*\eled)=\has^*(X,\eld).\ \ \hfill\Box$$

%De mani\`ere analogue, lorsque $A$ est un \fib \inve sur $X$, on peut construire un \mor
%\begin{equation}
%\label{ecuatia2000}
%can:\es^{n-2}\has^*(X,A)\otimes\esdoi\has^*(X,L\otimes A)\to\has^*(\hilx,\esdoi\elen\otimes\dna)
%\end{equation}
%
%
%\begin{teor}
%\label{teorema11}
%Soit $X$ une surface projective et lisse. Soit $L$ un fibr\'e inversible sur $X$. L'application canonique (\ref{ecuatia18})  pour $n=2$:
%\begin{equation}
%\label{ecuatia19}
%can:\esdoi\has^*(X,L)\to\has^*(\hild,\esdoi\eled)
%\end{equation}
%induit la d\'ecomposition en somme directe:
%$$\has^*(\hild,\esdoi\eled)=\esdoi\has^*(X,L)\oplus(\has^*(X,\maO_X)/\comp)\otimes\has^*(X,\eld).$$
%%%Elle est injective et son conoyau est 
%%%$$\bigoplus_{i\ge 1}\has^i(\maO_X)\otimes\has^{*-i}(\eld).$$
%
%\end{teor}
%
%
%\begin{cor}
%\label{cor12}
%
%Si $X$ est une vari\'et\'e projective lisse et de plus l'irr\'egularit\'e $q=0$ et $p_g=0$, alors l'application canonique (\ref{ecuatia19}) est un isomorphisme.
%
%\end{cor}

\vs Dor\'enavant, on omettera par convention l'espace $X$ dans la notation $ \has^*(X,L)$ pour tout faisceau inversible $L$.

\vs\par{\bf Preuve du \thm \ref{teorema11}:}

La proposition \ref{prop8} et les lemmes \ref{lema9} et \ref{lema10} impliquent la suite longue de cohomologie:
\begin{eqnarray}
\ldots\to\esdoi\has^*(L)=\has^*(\esxd,\dedoi)\to\has^*(\hild,\esdoi\eled)\to\esdoi\has^*(L)\oplus \has^*(\maO_X)\otimes\has^*(\eld)\to\has^*(\eld)\to\nonumber\\
\label{ecuatia20}
\to\has^{*+1}(\hild,\esdoi\eled)\to\ldots
\end{eqnarray}
La composition du morphisme (\ref{ecuatia002bis}) pour $k=1, n=2$ avec l'application 
$$\has^*(\hild,\esdoi\eled)\to\esdoi\has^*(L)$$
de la suite (\ref{ecuatia20}), provient de la composition 
$$\esdoi\has^*(L)=\has^*(\esxd,\dedoi)\stackrel{can}{\longrightarrow}\has^*(\hild,\esdoi\eled)\to\has^*(\hild,\dedoi).$$
Par la d\'efinition de ces morphismes il r\'esulte que cette composition co\"{\i}ncide avec l'isomorphisme  (\ref{ecuatia002bis}) pour $k=1$, $n=2$.
On obtient la d\'ecomposition en somme directe
$\has^*(\hild,\esdoi\eled)=\esdoi\has^*(L)\oplus K^*$, et la suite (\ref{ecuatia20}) se r\'eduit \`a 
\begin{equation}
\label{ecuatia20000}
\cdots\to K^*\to\has^*(\maO)\otimes\has^*(\eld)\to \has^*(\eld)\to K^{*+1}\to\cdots
\end{equation}
L'application 
$$\has^*(\maO_X)\otimes\has^*(\eld)\stackrel{b}{\to}\has^*(\eld)$$
de la suite (\ref{ecuatia20}) est celle canonique. Par cons\'equent cette application est surjective, donc $K^*=\Ker b$, d'o\`u la conclusion.
%% et la suite (\ref{ecuatia20}) se s\'epare en  suites exactes courtes:
%%\begin{eqnarray}\label{ecuatia21}
%%0\to\has^*(\hild,\esdoi\eled)\to\esdoi\has^*(L)\oplus \has^*(\maO_X)\otimes\has^*(\eld)\to\has^*(\eld)\to 0.
%%\end{eqnarray}
%%
%% Alors on peut r\'e-\'ecrire la suite exacte  (\ref{ecuatia21}) sous la forme:
%%\begin{equation}\label{ecuatia22}
%%0\to\esdoi\has^*(L)\to \has^*(\hild,\esdoi\eled)\to\has^*(\maO_X)\otimes\has^*(\eld)\to\has^*(\eld)\to 0.
%%\end{equation}
%%
%%
%%
%%On a vu que la derni\`ere application \'etait surjective. Son noyau est 
%%$$\bigoplus_{i\ge 1}\has^i(\maO_X)\otimes\has^{*-i}(\eld).$$
%%
%%Par cons\'equent la suite exacte (\ref{ecuatia22}) se r\'e-\'ecrit sous la forme \'enonc\'ee dans le \thm \ref{teorema11}. 
\hfill $\Box$

\vs\par{\bf Preuve du \th \ref{teor2.4'}:}

On tensorise la suite exacte (\ref{ecuatia17}) par le \fs \inve $\dda$ sur $\hild$. De la m\^eme mani\`ere que les \'enonc\'es \ref{prop8}, \ref{lema9}, \ref{lema10} on prouve
$$\has^*(\hild,\elded\otimes\dda)=\has^*(X,A)\otimes\has^*(X,\eld\otimes A),$$
$$\has^*(\hild,\dedoi\otimes\dda)=\esdoi\has^*(X,L\otimes A),$$ 
$$\has^*(\hild,\mu^*\eld\vert_{\dhild}\otimes\dda)=\has^*(X,\eld\otimes A^2).$$
On obtient la suite exacte longue:
\begin{eqnarray}
\cdots\to \has^*(\hild,\esdoi\eled\otimes\dda)&\to& \esdoi\has^*(X,L\otimes A)\oplus \has^*(X,A)\otimes\has^*(X,\eld\otimes A)\to\nonumber\\
\label{ecuatia2002}
&\to&\has^*(X,\eld\otimes A^2)\to\cdots
\end{eqnarray}
De m\^eme que dans la preuve du \th \ref{teorema11}, la composition
$$\esdoi\has^*(X,L\otimes A)\stackrel{can}{\to}\has^*(\hild,\esdoi\eled\otimes\dda)\to \esdoi\has^*(X,L\otimes A)$$
est un isomorphisme. Alors $\has^*(\hild,\esdoi\eled\otimes\dda)=\esdoi\has^*(X,L\otimes A)\oplus K^*$, et la suite exacte (\ref{ecuatia2002}) se transforme dans la suite exacte (\ref{ecuatia2001}). \hfill $\Box$

\section{ Le cas $n=3$, le passage de $\hilts$ \`a $\hilt$}
\label{sectiune3}

%On note encore $\mu:\dhilt\to X$ l'application qui associe au sch\'ema $Z$ de longueur $3$ son unique point multiple. Elle prolonge, dans le cas $n=3$, l'application $\mu:\dhilts\to X$ d\'efinie dans la section \ref{sectiune1}.
%Le but de cette section est de d\'emontrer:

Le but de cette section est de prolonger la suite exacte  (\ref{ecuatia000}) de la vari\'et\'e $\hilts$ \`a une suite exacte de faisceaux sur $\hilt$. Le r\'esultat est r\'esum\'e dans le \th suivant. Avant de l'\'enoncer on introduit quelques notations.

La sous-vari\'et\'e des sch\'emas singuliers $\dhilt$ est une hypersurface de la vari\'et\'e lisse $\hilt$. L'ouvert $\dhilts=\dhilt\cap\hilts$ est lisse et irr\'eductible, et son compl\'ementaire, $\dhilt\setminus\hilts$, est de codimension $2$ dans $\hilt$. Par cons\'equent $\dhilt$ est une vari\'et\'e r\'eduite et irr\'eductible, donc int\`egre. On consid\`ere sa normalisation $\pi:\Sigma\to\dhilt$. L'ouvert $\dhilts$ est lisse, donc $\pi$ est un \iso au-dessus de $\dhilts$. On note $i:\dhilts\to\Sigma$ le \mor d'inclusion.

\begin{lema}
\label{lema300}
On consid\`ere le diagramme:
\begin{equation}
\label{ecuatia301}
{\diagram
&\Sigma\dlto_{\pi}&\\
\dhilt&\dhilts\lto_{j}\uto_{i}\rto^{\mu}&X,
\enddiagram}
\end{equation}
o\`u $\mu$ est le \mor d\'efini dans la section \ref{sectiune1}. Il existe une application r\'eguli\`ere $\tmu:\Sigma\to X$ qui rend commutatif ce diagramme.
\end{lema}

La preuve du lemme \ref{lema300} sera donn\'ee apr\`es la preuve du \th \ref{teorema13}.

\begin{teor}
\label{teorema13}

La suite exacte (\ref{ecuatia000}) de faisceaux sur $\hilts$  se prolonge \`a une suite exacte sur $\hilt$:
$$0\to\esdoi\eltr\to\dedoit\oplus\eldtr\to\pi_*\tmu^*\eld\to 0.$$

\end{teor}

\vs\par {\bf Preuve du \thm \ref{teorema13}:}

La preuve utilise les r\'esultats suivants, dont la d\'emonstration occupera le reste de la section:

%l'extension du diagramme (\ref{diag16}) de $\hilts$ \`a $\hilt$.
%On commence avec la proposition technique:

\begin{prop}
\label{prop14}

Le faisceau $\dedoit$ est localement libre de rang $3$ sur $\hilt$.

\end{prop}

\begin{lema}
\label{lema19}

Le morphisme 
$\eldtr\to\mu^*\eld\vert_{\dhilts}$
sur $\hilts$   de la suite exacte (\ref{ecuatia000}) se prolonge \`a un morphisme sur $\hilt$  
$$\eldtr\to\pi_*\tmu^*\eld.$$

\end{lema} 

\vs On consid\`ere le ferm\'e $X\subset\hilt$ des sch\'emas d'id\'eal $m_x^2, x\in X$.

\begin{lema}
\label{lema20}

Le morphisme  
$\dedoit\to\mu^*\eld\vert_{\dhilts}$
sur $\hilts$ de la suite exacte (\ref{ecuatia000}) se prolonge \`a un morphisme \begin{equation}
\label{ecuatia301,7}
\dedoit\to\pi_*\tmu^*\eld
\end{equation}
 sur $\hilt$, surjectif sur $\hilt\setminus X$.

\end{lema} 

\vs On note $k:\hilts\to\hilt\setminus X$ l'inclusion canonique. Le foncteur $k_*$ est exacte \`a gauche. On l'applique \`a la suite exacte (\ref{ecuatia000}). On obtient, en tenant compte du fait que les faisceaux $\esdoi\eltr$, $\eldtr$ et $\dedoit$ sont localement libres (on utilise ici la proposition \ref{prop14}), la suite exacte sur $\hilt$:
\begin{equation}
\label{ecuatia302}
0\to\esdoi\eltr\to\dedoit\oplus\eldtr\stackrel{b}{\to} k_*\mu^*\eld.
\end{equation}
Le \mor $i$ est une immersion ouverte, donc le diagramme (\ref{ecuatia301}) induit une inclusion canonique $\tmu^*\eld\subset i_*\mu^*\eld.$ On lui applique le \mor fini $\pi$ et on obtient l'inclusion canonique
$$\pi_*\tmu^*\eld\subset \pi_*i_*\mu^*\eld=k_*\mu^*\eld.$$
Les lemmes \ref{lema19} et \ref{lema20} montrent que le \mor $b$ de la suite (\ref{ecuatia302}) se factorise \`a travers un \mor
$$\dedoit\oplus\eldtr\to\pi_*\tmu^*\eld\subset k_*\mu^*\eld.$$
La suite exacte (\ref{ecuatia302}) induit par cons\'equent la suite exacte sur $\hilt\setminus X$:
$$0\to\esdoi\eltr\to\dedoit\oplus\eldtr\to\pi_*\tmu^*\eld.$$
L'application $b$ de la suite exacte (\ref{ecuatia302}) est surjective d'apr\`es le lemme \ref{lema20}.
On note $j:\hilt\setminus X\to\hilt$ l'inclusion canonique. On applique le foncteur $j_*$ \`a la derni\`ere suite exacte et on trouve une suite exacte sur $\hilt$:
\begin{equation}
\label{ecuatia850}
0\to\esdoi\eltr\to\dedoit\oplus\eldtr\to j_*\pi_*\tmu^*\eld\to \er^1j_*(\esdoi\eltr)
\end{equation}
Le \mor $j$ est l'inclusion d'un ouvert dont le compl\'ementaire est de codimension $4$ dans la vari\'et\'e lisse $\hilt$, et $\esdoi\eltr$ est un \fs localement libre sur $\hilt$, donc $\er^1j_*(\esdoi\eltr)=0$. Pour \'ecrire la suite (\ref{ecuatia850}) sous la forme de l'\'enonc\'e il suffit de d\'emontrer
\begin{equation}
\label{ecuatia309}
\pi_*\tmu^*\eld=j_*j^*(\pi_*\tmu^*\eld).
\end{equation}
On consid\`ere le diagramme
$$\diagram
\Sigma'\rto^{l}\dto_{\pi}&\Sigma\\
\hilt\setminus X\rto^{j}&\dhilt,
\enddiagram$$
o\`u $\Sigma'$ est l'image r\'eciproque de $\hilt\setminus X$ par le \mor $\pi$. On a
$$j_*j^*\pi_*\tmu^*\eld=j_*\pi_*l^*\tmu^*\eld=\pi_*l_*l^*\tmu^*\eld.$$
Le faisceau $\tmu^*\eld$ est inversible sur $\Sigma$ et $l$ est l'inclusion d'un ouvert dont le compl\'ementaire est de codimension $3$ dans la vari\'et\'e normale $\Sigma$. On en d\'eduit
$$l_*l^*\tmu^*\eld=\tmu^*\eld,$$
d'o\`u l'\'egalit\'e (\ref{ecuatia309}), d'o\`u la suite exacte de l'\'enonc\'e. \hfill $\Box$
%
%
%
%La surjectivit\'e du \mor (\ref{ecuatia301,7}) termine la preuve du \th \ref{teorema13}. \hfill $\Box$

\vs Le reste du chapitre est consacr\'e \`a la d\'emonstration du lemme \ref{lema300}, de la proposition \ref{prop14} et des lemmes \ref{lema19}, \ref{lema20} dans cet ordre.

\vs\par {\bf Preuve du lemme \ref{lema300}:}

On utilisera une pr\'esentation explicite de la normalisation $\pi:\Sigma\to\dhilt.$ On consid\`ere l'image r\'eciproque $D$ de l'hypersurface $\dhild$ par le \mor $p_{32}:\hiltd\to\hild.$ C'est le ferm\'e des points $(Z,Z')\in\hilt\times\hild$ qui v\'erifient $Z'\subset Z$ et $Z'$ est singulier. Le \mor $\pi_{32}:\hiltd\to\hilt$ envoie $D$ dans $\dhilt$ et ce \mor est un \iso au-dessus de $\dhilts$. Effectivement, un point de $\dhilts$ est de la forme $\maO_Z=\maO_{Z'}\oplus\maO_a$, avec $Z'$ sch\'ema singulier de longueur $2$ \`a support disjoint du point $a$. L'association $Z\mapsto(Z,Z')$ est une application r\'eciproque sur $\dhilts$ de l'application $\pi_{32}:D\to\dhilt.$

Consid\'erons la factorisation de Stein du \mor $\pi_{32}:D\to\dhilt$:
$$D\stackrel{\tpi_{32}}{\longrightarrow}\Sigma\stackrel{notation}{=}\Spec(\pi_{32*}\maO_{D})\stackrel{\pi}{\to}\dhilt.$$
On d\'emontre dans la suite que le sch\'ema $D$ est normal. Le fait que pour tout ouvert $U$ de $\Sigma$ on a $\Gamma(U,\maO_{\Sigma})=\Gamma(\tpi^{-1}_{32}(U),\maO_D)$ implique que $\Sigma$ est normal. Puisque le \mor $\pi_{32}$ est birationnel, le \mor $\pi$ est aussi birationnel. Donc le \mor $\pi:\Sigma\to\dhilt$  est un \mor de normalisation.

On construit le diagramme:
\begin{equation}
\label{ecuatia302,5}
{\diagram
D\dto_{\tpi_{32}}\rto^{p_{32}}\drto^{\ttmu}&\dhild\dto^{HC}\\
\Sigma&X=\desxd,
\enddiagram}
\end{equation}
o\`u $HC$ est le \mor de Hilbert-Chow, et $\ttmu$ est la composition $HC\circ p_{32}$. Le \mor $\ttmu$ associe au couple $(Z,Z')$ le point singulier du sch\'ema $Z$. Il se factorise alors \`a travers une application continue $\dhilt\to X$, donc \`a travers  une application continue $\Sigma\stackrel{\tmu}{\to}X$. Cette derni\`ere application est r\'eguli\`ere, car pour tout ouvert $W\subset X$ on dispose d'un \mor d'anneaux:
$$\Gamma(W,\maO_X)\to\Gamma(\tmu^{-1}(W),\maO_{\Sigma})=\Gamma(\tpi^{-1}_{32}\tmu^{-1}(W),\maO_{D}).$$ 
En consid\'erant les isomorphismes $D\stackrel{\tpi_{32}}{\to}\Sigma\stackrel{\pi}{\to}\dhilt$ au-dessus de $\dhilts$, on voit que le \mor qu'on vient de construire prolonge le \mor $\mu:\dhilts\to X$ construit dans la section \ref{sectiune1}. Pour d\'emontrer que le sch\'ema $D$ est normal, on consid\`ere le \mor $\pi_{32}:D\to\dhilt.$ On prend le ferm\'e $X\subset\dhilt$ des sch\'emas d'id\'eal $m^2_x$, $x\in X$. La fibre du \mor $\pi_{32}:D\to\dhilt$ dans chaque point $\{m^2_x\}$ est de dimension $1$: elle est form\'ee par les points $Z'\in\hild$, sch\'emas de longueur $2$ d'id\'eal $\comp\cdot \ell+m^2_x, \ell\ne 0$. Par suite l'image r\'eciproque $Y=\pi_{32}^{-1}(X)$ est de dimension $3$.

On d\'emontre:
\begin{lema}
\label{lema310}
La vari\'et\'e $D\setminus Y$ est lisse.
\end{lema}

La vari\'et\'e $\hiltd$ est lisse d'apr\`es \cite{Che}, \cite{Tik} (la preuve  du lemme \ref{lema16}(ii) qui suit donne une autre d\'emonstration), et $D$ est une hypersurface dans $\hiltd$. D'apr\`es le lemme \ref{lema310} la vari\'et\'e $D$ est r\'eguli\`ere en codimension $1$. Par cons\'equent (\cite{Hart}, prop. II 8.23.b) la vari\'et\'e $D$ est normale. Ceci termine la preuve du lemme \ref{lema300}. \hfill $\Box$

\vs\par {\bf Preuve du lemme \ref{lema310}:}
  
La question est locale et il suffit de la traiter en g\'eom\'etrie analytique. On peut supposer que $X=\pp$. On note $(U:V:W)$ les coordonn\'ees homog\`enes sur $\pp$ et $u=\frac{U}{W}, v=\frac{V}{W}$. Le \mor $\pi_{32}:D\to\dhilt$ est un \iso au-dessus de $\dhilts$, et $\dhilts$ est lisse. Il suffit donc de traiter le probl\`eme au-dessus de $\dhilt\setminus(X\cup\dhilts).$ Le groupe $\pgl(3)$ agit sur cet ensemble ($X=\pp$) et il y a deux orbites, repr\'esent\'ees par les sch\'emas d'id\'eal $(v,u^3)$, respectivement $(v-A u^2, u^3), A\in\comp^*$.

\begin{lema}
\label{lema3.600}
L'association:
\begin{equation}
\label{ecuatia303}
(A,B,C,p,q,r)\mapsto(v-(Au^2+Bu+C),u^3+pu^2+qu+r)
\end{equation}
donne une param\'etrisation locale $\comp^6\to\hilt$ de $\hilt$ au voisinage du point $s_0=(v,u^3)$.
\end{lema}

\vs\par {\bf Preuve:}

L'application 
$$(A,B,C)\mapsto v-(Au^2+Bu+C)\vert_{\has^0(\maO_{s_0})}$$
est une application affine de rang $3$ entre les espaces vectoriels $\comp^3\to\comp^3$. Par suite l'application $(A,B,C)\mapsto
v-(Au^2+Bu+C)\vert_{\has^0(\maO_{s_0})}$ est de rang $3$ pour $s$ dans un voisinage de $s_0$ dans $\hilt$.
On note $\Gamma\simeq{\mathbb{A}}^3$ l'ensemble des coniques d'\'equation $v=Au^2+Bu+C$. Alors tout point $s$ dans un voisiange $V$ de $s_0$ dans $\hilt$ est contenu  une unique conique $\gamma(s)\in\Gamma$  et l'application $\gamma:V\to\Gamma$ ainsi obtenue admet des diff\'erentielles de rang $3$ en tout point. Les fibres du \mor $\gamma$ en chaque point $(A,B,C)\in\Gamma$ sont en bijection avec les sch\'emas de longueur $3$ sur la droite $v=0$. 
Ces sch\'emas sont param\'etr\'es au voisinage du point $u=0$ par $(p,q,r)\in\comp^3$ par l'association
$$(p,q,r)\mapsto \comp[u]/(u^3+pu^2+qu+r).$$
On obtient une application locale $V\to\comp^6, \ \ \ \ s\mapsto(A,B,C,p,q,r)$. La diff\'erentielle de cette application est de rang $6$ en $s_0$, donc il s'agit d'un \iso local. \hfill $\Box$

\vs Tout voisinage du point $s_0=(v,u^3)$ contient des points de la forme $(v-A u^2,u^3)$: on prend $B=C=p=q=r=0$. Il suffit donc de d\'emontrer que le sch\'ema $D$ est lisse au-dessus d'un voisinage du point $s_0=(v,u^3)$.

La vari\'et\'e $\dhild$ admet une param\'etrisation locale $\comp^3\to\dhild$ au voisinage du point $(v,u^2)$ donn\'ee par
\begin{equation}
\label{ecuatia304}
(\lambda,\mu,\alpha)\mapsto(v-\lambda u-\mu,(u-\alpha)^2).
\end{equation}
Un calcul simple prouve que le sch\'ema d'\'equations (\ref{ecuatia303}) contient le sch\'ema d'\'equations (\ref{ecuatia304}) si et seulement si:
%\begin{eqnarray}
%\label{ecuatia305}
%v-(Au^2+Bu+C)&=&P\cdot(v-\lambda u-\mu)+Q\cdot(u-\alpha)^2,\\
%\label{ecuatia306}
%u^3+pu^2+qu+r&=&R\cdot(v-\lambda u-\mu)+S\cdot(u-\alpha)^2,
%\end{eqnarray}
%pour $P,Q,R,S\in\comp[[u,v]]$. L'\'egalit\'e (\ref{ecuatia305}) \'equivaut \`a
%$P=1+T\cdot(u-\alpha)^2$, $ Q=-A-T\cdot(v-\lambda u-\mu)$, $T\in\comp[[u,v]]$ et
\begin{eqnarray}
-B=-\lambda-2A\alpha \nonumber\\
\label{ecuatia3078}
-C=-\mu-A\alpha^2\\
q=-2p\alpha-3\alpha^3\nonumber\\
r=-q\alpha-p\alpha^2-\alpha^3.\nonumber
\end{eqnarray}
%L'\'egalit\'e (\ref{ecuatia306}) \'equivaut \`a
%$R=F\cdot(u-\alpha)^2$, $S=(u+p+2\alpha)-F\cdot(\mu-\lambda u-\mu)$, $F\in\comp[[u,v]]$ et
%\begin{eqnarray}
%q=-2p\alpha-3\alpha^3, \nonumber \\
%\label{ecuatia308}
%r=-q\alpha-p\alpha^2-\alpha^3.
%\end{eqnarray}
En conclusion, le sous-sch\'ema $D$ est d\'efini par les \'equations (\ref{ecuatia3078}) dans $\hilt\times\dhild$
%les \'equations  (\ref{ecuatia307}) et (\ref{ecuatia308}) repr\'esentent les \'equations du ferm\'e $D$ dans $\hilt\times\dhild$ 
au voisinage du point $((v,u^2),(v,u^3))$. C'est le graphe d'une application $(A,p,\lambda,\mu,\alpha)\mapsto(B,C,q,r)$, donc le sous-sch\'ema $D$ est lisse dans ce voisinage. \hfill $\Box$

\vs\par{\bf Preuve de la proposition \ref{prop14}:}

La relation (\ref{ecuatia7}) du lemme \ref{lema5}:
$$\dedoit=\dedoi\boxtimes\maO\oplus \eled\boxtimes L_3$$
prouve que $\dedoit$ est localement libre de rang $3$ sur $\hilts$. Le ferm\'e de dimension $4$, $\hilt\setminus\hilts$ est la r\'eunion des ensembles disjoints $X\subset\hilt$ (qui param\`etre les points triples de la forme $\maO/_{\emfrac_X^2}$, avec $x\in X$), et $\hiltc$ (qui param\`etre les sous-sch\'emas de longueur $3$ \`a support en un point et situ\'es sur une courbe lisse dans $X$).
Au-dessus de $\hiltc$ la projection $\pi_{32}$ qui d\'efinit $\dedoit$ (cette projection a \'et\'e introduite dans la section \ref{sectiune1}), 
est quasi-finie. Alors la dimension de la fibre de $\dedoit$ en un point de $\hiltc$ co\"{\i}ncide avec la longueur de la fibre sch\'ematique de $\pi_{32}$ en ce point. La longueur de la fibre sch\'ematique d\'epend seulement de la g\'eom\'etrie analytique de $X$. C'est par cons\'equent une constante le long de $\hiltc$.
De m\^eme, la dimension de la fibre de $\dedoit$ en un point $x\in X\subset\hilt$ ne d\'epend pas du choix de $x\in X$. Pour d\'emontrer la proposition il suffit de d\'emontrer que la dimension de la fibre de $\dedoit$ en un point de $x\in X$ est $3$. Effectivement, tout voisinage de $x$ dans $\hilt$ contient des points de $\hiltc$ et par le \thm de semi-continuit\'e (\cite{Hart}, ex. III 12.7.2) on obtient que la dimension des fibres de $\dedoit$ aux points de $\hiltc$ est $3$.

Il suffit donc de d\'emontrer que $\dedoit\otimes_{\maO_{\hilt}}\maO_X$ est localement libre de rang $3$ le long de $X\subset\hilt$.
On d\'efinit le ferm\'e $Y$ dans $\hiltd$ par le diagramme cart\'esien
$$\diagram
Y\rto^{j}\dto_{\pi'_{32}}&\hiltd\dto_{\pi_{32}}\\
X\rto^{i}&\hilt.
\enddiagram$$

La formule de changement de base dans la cat\'egorie d\'eriv\'ee pour l'application $i:X\to\hilt$ et pour le faisceau $\eldoi=p^*_{32}\dedoi$ sur $\hiltd$:
$$\er\pi_{32*}(\eldoi)\otimes^L_{\hilt}\maO_X=\er\pi'_{32*}(\eldoi\otimes^L_{\hilt}\maO_X)$$
induit les suites spectrales:
\begin{eqnarray*}
\grave{ }E_2^{pq}=\er^p\pi'_{32*}(\underline{\Tor}^{\maO_{\hilt}}_{-q}(\maO_X,\eldoi))\Rightarrow\te^{p+q}\eldoi\\
''E_2^{pq}=\underline{\Tor}^{\maO_{\hilt}}_{-p}(\maO_X,\er^q\pi_{32*}\eldoi)\Rightarrow\te^{p+q}\eldoi.
\end{eqnarray*}

On commence par calculer les termes $\grave{ }E_2^{pq}$. Le faisceau $\eldoi$ est inversible sur $\hiltd$. Alors 
$$\underline{\Tor}^{\maO_{\hilt}}_{-q}(\maO_X,\eldoi)=\underline{\Tor}^{\maO_{\hilt}}_{-q}(\maO_X,\maO_{\hiltd})\otimes\eldoi\vert_Y.$$
\def\pipstd{\pi^{\prime *}_{32}}
\begin{lema}
\label{lema15}

Le faisceau $\eldoi\vert_Y$ est isomorphe \`a $\pipstd\eld$.

\end{lema}

\vs\par{\bf Preuve:}

Le diagramme
$$\diagram
Y\rto^{j}\drto_{\pi'_{32}}&\hiltd\rto^{p_{32}}&\hild\rto&\esxd\\
&X\urrto_{diag}&&
\enddiagram$$
est commutatif. Alors $j^*\eldoi=\pipstd diag^*\dedoi=\pipstd\eld.\ \ \hfill\Box$

Pour calculer $\grave{ }E_2^{pq}$ on utilise le lemme technique:

\begin{lema}
\label{lema16}
\begin{itemize}
\item{\rm (i)} Le ferm\'e $Y$ est une fibration en espaces projectifs $\proj_1$ au-dessus de $X$.
\item{\rm (ii)} Le faisceau ${\uTor}^{\maO_{\hilt}}_{q}(\maO_X,\maO_{\hiltd})$ est nul pour $q>1$, $\maO_Y$ pour $q=0$ et un faisceau inversible $M$ sur $Y$ pour $q=1$.
\item{\rm (iii)} Au-dessus de  chacune des fibres $\proj_1$ du morphisme $\pi'_{32}:Y\to X$, $M$ est isomorphe \`a $\maO(-3)$.
\end{itemize}
\end{lema}

La d\'emonstration du lemme \ref{lema16} sera donn\'ee \`a la fin de la d\'emonstration de la proposition \ref{prop14}.

On a 
$$\begin{array}{ccccc}
   \grave{ }E_2^{p0}&=&\er^p\pi'_{32*}(\pipstd\eld)&=&\hfill\\
\hfill&=&\eld\otimes\er^p\pi'_{32*}\maO_Y&=&
\begin{cases}
  0 \mbox{ \ \ \ si \ \ \ } p>0\cr
\eld \mbox{ \ \ \ si \ \ \ } p=0,
\end{cases}
 \end{array}$$
d'apr\`es le point (i) du lemme.
On a aussi
$$\begin{array}{ccccc}
   \grave{ }E_2^{p,-1}&=&\er^p\pi'_{32*}(\pipstd\eld\otimes M)&=&\hfill\\
\hfill&=&\eld\otimes\er^p\pi'_{32*}M&=&
\begin{cases}
  0 \mbox{ \ \ \ si \ \ \ } p\ne 1\cr
 \mbox{ localement\ \ libre \ \ de \ \ rang \ \ 2\ \ sur\ \ X\ \ \ si \ \ \ } p=1,
\end{cases}
 \end{array}$$
d'apr\`es l'\'egalit\'e du point (iii) du lemme et du fait que $h^1(\proj_1,\maO(-3))=2.$
Finalement, le point (ii) du lemme implique 
$$\grave{ }E_2^{pq}=0 \ \ \mbox{\ \ \ pour\ \ \ }q\ne 0,-1\ .$$
En conclusion, on a $\te^{p+q}\eldoi=0$ pour $p+q\ne 0$ et $\te^0\eldoi$ est un faisceau localement libre de rang $3$  sur $X$ qui v\'erifie la suite exacte:
\begin{equation}
\label{ecuatia23}
0\to\eld\otimes\er^1\pi'_{32*} M\to\te^0\eldoi\to\eld\to 0.
\end{equation}

On consid\`ere maintenant $''E^{pq}_2$. On d\'emontre tout d'abord:
\begin{equation}
\label{ecuatia24}
\er^q\pi_{32*} \eldoi=0 \ \ \mbox{\ \ \ pour\ \ \ }q>0\ .
\end{equation}

On suppose par l'absurde qu'il existe un entier $q>0$ tel que le faisceau $\er^q\pi_{32*}\eldoi$ est non nul. On consid\`ere le plus grand parmi ces entiers $q$. Le faisceau $\er^q\pi_{32*}\eldoi$ a son support sur $X\subset\hilt$ puisque au dessus de $\hilt\setminus X$ le morphisme $\pi_{32}$ est fini. 
Alors
$$''E_2^{0q}=\maO_X\otimes_{\maO_{\hilt}}\er^q\pi_{32*}\eldoi\ne 0.$$
Mais le fait que $q$ est le plus grand entier tel que  $\er^q\pi_{32*}\eldoi\ne 0$ implique que les diff\'erentielles $''d_2,''d_3,\ldots$  qui aboutissent sur $''E_2^{0q}$ et qui partent de $''E_2^{0q}$ sont nulles. 
Donc $\te^q\eldoi=''E_2^{0q}\ne 0$, contradiction.
Par cons\'equent les seuls termes non nuls sont $''E_2^{p0}$. Donc $''E_2^{p0}=\te^p\eldoi$. Alors le seul terme non nul est 
$$\te^0\eldoi=''E^{00}_2=\maO_X\otimes_{\maO_{\hilt}}\pi_{32*}\eldoi=\maO_X\otimes_{\maO_{\hilt}}\dedoit.$$

Alors le faisceau $\dedoit\otimes_{\hilt}\maO_X$ est localement libre de rang $3$ sur $X$ et v\'erifie la suite exacte:
\begin{equation}
\label{ecuatia25}
0\to\eld\otimes\er^1\pi'_{32}M\to\dedoit\otimes_{\hilt}\maO_X\to\eld\to 0.
\end{equation}
Cela termine la d\'emonstration de la proposition \ref{prop14}. \hfill $\Box$

\vs\par {\bf Preuve du lemme \ref{lema16}:}

%\begin{itemize}
%\item
 {\rm (i)} Soit $x\in X\subset\hilt$. Les points de la fibre en $x$ du morphisme $\pi_{32}:\hiltd\to\hilt$ sont donn\'es par les sch\'emas de longueur $2$ inclus dans $s=\Spec \maO_X/_{\emfrac^2_x}$. Ces sous-sch\'emas sont de la forme $b=\Spec \maO_X/_{\emfrac^2_x,\ell_x}$, pour une droite $\ell_x$ qui passe par $x$. On d\'eduit alors que g\'eom\'etriquement la fibre en $x$ du morphisme $\pi_{32}$ est $\proj_1$. Pour d\'emontrer que c'est aussi sch\'ematiquement $\proj_1$ il suffit de d\'emontrer que cette fibre est lisse en tout point.
L'\'enonc\'e est local et analytique. On peut supposer que $X=\proj_2$, avec les coordonn\'ees $(U:V:W)$,  $x$ est le point $(0:0:1)$ et $\ell_x=V$, donc $b=\Spec \maO_{\proj_2}/_{(U^2,V)}$. On utilise l'affirmation suivante, tautologique. Soit $L$ un fibr\'e inversible sur $X$ et $s\in\has^0(X,L)$ une section globale.   
On consid\`ere le diagramme
$$\diagram
\Xi\rto^{p_{21}}\dto_{\pi_{21}}&X\\
\hild&.
\enddiagram$$
Par d\'efinition $\eled=\pi_{21*}p_{21}^*L$.  L'image r\'eciproque $p^*_{21}s\in\has^0(\Xi,p_{21}^*L)$ est une section globale sur $\Xi$ et s'identifie \`a une section $[s]\in\has^0(\hild,\eled)$.

\begin{rem}
\label{obs17}
{\rm
Soit $z\in\hild$ un point qui correspond \`a un sous-sch\'ema de longueur $2$, $Z\subset X$. Le sous-sch\'ema $Z\subset X$ est contenu dans la courbe d'\'equation $s=0$ si et seulement si la section $[s]\in\has^0(\hild,\eled)$ s'annule au point $z\in\hild$.}

\end{rem}

\vs\par{\bf Preuve:}

La section $[s]$ s'annule au point $z$ si et seulement si $p^*_{21}s$ s'annule sur la fibre ${Z}\subset\Xi$. Ceci \'equivaut \`a l'annulation de $s$ en ${Z}\subset X$. \hfill $\Box$

\vs Ici $L=\maO(2)$. Par suite, comme $s$ est d'id\'eal $(U^2,UV, V^2)$ dans $\pp$, la fibre en $s$ du morphisme $\pi_{32}$ est donn\'ee par le ferm\'e d'\'equations $[U^2],[UV],[V^2]$ dans $\ppd$.
Le point $b$ admet la param\'etrisation locale suivante:
\begin{eqnarray}
\label{ecuatia26}
\comp^4&\to&\ppd\nonumber\\
(\alpha,\beta,\lambda,\mu)&\to&\maI_Z=(V-\lambda U-\mu W, U^2-\alpha UW-\beta W^2).
\end{eqnarray}
Les sections $[UW]$ et $[W^2]$ engendrent le \fs localement libre $\eled$ au voisinage de $b$ et on a:
\begin{eqnarray*}
{[U^2]}&=&\alpha[UW]+\beta[W^2]\\
{[UV]}&=&(\lambda\alpha+\mu)[UW]+\lambda\beta[W^2]\\
{[V^2]}&=&(\lambda^2\alpha+2\lambda\mu)[UW]+(\lambda^2\beta+\mu^2)[W^2].
\end{eqnarray*}

On obtient que la fibre en $s$ est donn\'ee dans la param\'etrisation de $b$ par les \'equations:
$$\alpha=\beta=\lambda\alpha+\mu=\lambda\beta=\lambda^2\alpha+2\lambda\mu=\lambda^2\beta+\mu^2=0,$$
\'equivalentes \`a 
$$\alpha=\beta=\mu=0.$$
Par cons\'equent la fibre en $s$ est lisse de dimension $1$ au voisinage de $b$, ce qu'il fallait d\'emontrer.

%\item 
 {\rm (ii)} On applique le r\'esultat connu:
\begin{lema}
\label{lema18}

Soit $B$ une alg\`ebre, $A$ une $B$-alg\`ebre plate, $M$ un $A$-module et $N$ un $B$-module. Alors:
$$\Tor^A_q(A\tens_B N,M)=\Tor^B_q(N,M),$$

\end{lema}

pour $A=\maO_{\hilt\times\hild}, B=\maO_{\hilt}, M=\maO_{\hiltd}, N=\maO_X$. On trouve:
$$\uTor^{\maO_{\hilt}}_q(\maO_{\hiltd},\maO_X)=\uTor^{\maO_{\hilt\times\hild}}_q(\maO_{\hiltd},\maO_{X\times \hild}).$$
On calculera le membre de droite en \'ecrivant explicitement les \'equations de $\hiltd$  et $X\times\hild$ dans $\hilt\times\hild$. On se place au voisinage du point $ (s,b) \in \hilt\times\hild$ consid\'er\'e au point (i).
Le point $b$ admet la param\'etrisation locale (\ref{ecuatia26}). 
Le lemme suivant donne une param\'etrisation locale du point $s\in\hilt$:

\begin{lema}
\label{lema3.11}
Consid\'erons une application lin\'eaire $m:\comp^2\to\comp^3$. Les mineurs $2\times 2$ de l'application lin\'eaire 
\begin{equation}
\label{ecuatia2999}
a(m)=\left(\begin{array}{cc}V&0\\ -U&V\\ 0&-U\end{array}\right)+\left(\begin{array}{cc}m_1&m_4\\ m_2&m_5\\ m_3&m_6\end{array}\right)W
\end{equation}
d\'efinissent l'id\'eal d'un point $s_m$ de $\hilt$. Le \mor $m\mapsto s_m$ 
est un \iso d'un voisinage de $0$ dans $L(\comp^2,\comp^3)$ sur un voisinage de $s$ dans $\hilt$.
\end{lema}

\vs\par{\bf Preuve:}

Prenons le complexe
\begin{equation}
\label{ecuatia3000}
0\to\maO(-3)^2\stackrel{a}{\to}\maO(-2)^3\stackrel{b}{\to}\maO\to 0,
\end{equation}
o\`u $b$ est d\'efini par les mineurs $2\times 2$ de l'application $a$. Le \mor $a$ est injectif et par construction on a $\Ker b=\Im a$. Le faisceau $\maO/\Im b$ est de rang $0$, $c_1=0$ et $\chi=3$. Il est par suite l'anneau structural d'un sch\'ema $s_m$ de longueur $3$ sur $X$. On obtient alors un \mor $m\mapsto s_m$ de $L(\comp^2,\comp^3)$ dans $\hilt$. Pour la construction du \mor r\'eciproque on consid\`ere l'ouvert $\hilt\setminus(\maH\cup\partial_W)$, o\`u $\maH$ est l'hypersurface des sch\'emas dont le support se trouve sur une droite et $\partial_W$ est le ferm\'e des sch\'emas dont le support rencontre la droite \`a l'infini $W=0$. D'apr\`es la suite spectrale de Beilinson (\cite[page 240]{O-S-S} ) on trouve que l'id\'eal de tout point de $\hilt\setminus\maH$ est le conoyau d'un \mor injectif $a:\maO(-3)^2\to\maO(-2)^3$, donn\'e par une matrice $a=A+B\cdot W$, o\`u $A$ est une forme lin\'eaire en $U$ et $V$, et $B$ une matrice constante. Si en plus le point se trouve dans $\hilt\setminus\partial_W$, le \mor $a$ doit \^etre injectif en restriction \`a chaque point de la droite $W=0$. Par cons\'equent l'application lin\'eaire d\'efinie par la matrice $A$ est injective en restriction  \`a tout point de la droite $W=0$. Alors la matrice $A$ est conjugu\'ee, par un unique \'el\'ement du groupe $\gl(3)\times\gl(2)/\comp^*$, avec la matrice 
$$A_0=\left(\begin{array}{cc}V&0\\ -U&V\\ 0&-U \end{array}\right).$$
En conclusion, l'id\'eal de tout point de  $\hilt\setminus(\maH\cup\partial_W)$ est le conoyau d'un \mor injectif $a=A_0+B_0\cdot W, B_0\in L(\comp^2,\comp^3)$. L'application $\hilt\setminus(\maH\cup\partial_W)\to L(\comp^2,\comp^3)$ donn\'ee par $z\mapsto B_0$ constitue le \mor r\'eciproque recherch\'e. \hfill $\Box$

\vs Donc  le point $s\in\hilt$ admet la param\'etrisation locale suivante:
$$(m_1,m_2,m_3,m_4,m_5,m_6)\in\comp^6\to\maI_Z=(M_{23},M_{13},M_{12}),$$
o\`u $M_{23}, M_{13}, M_{12}$ sont les mineurs $2\times 2$ de la matrice (\ref{ecuatia2999}).
%%%
%%%
%%%\begin{eqnarray}
%%%\label{ecuatia27}
%%%M_1=\left|\begin{array}{cc}-U+m_2W&m_3W\\V+m_5W&-U+m_6W\end{array}\right|, M_2=\left|\begin{array}{cc}V+m_1W&m_3W\\m_4W&-U+m_6W\end{array}\right|, \nonumber \\
%%%M_3=\left|\begin{array}{cc}V+m_1W&-U+m_2W\\m_4W&V+m_5W\end{array}\right|
%%%\end{eqnarray}
%%%sont les mineurs $2\times 2$ de la matrice
%%%$$\left(\begin{array}{ccc}V&-U&0\\ 0&V&-U\end{array}\right)+\left(\begin{array}{ccc}m_1&m_2&m_3\\ m_4&m_5&m_6\end{array}\right)W.$$
Au voisinage de $s$, les points de $\hilt$ de la forme $\emfrac_{y}^2$ ont pour id\'eal $(U-m_2W,V+m_1W)^2$ et sont donn\'es par les \'equations:
\begin{equation}
\label{ecuatia28}
m_1-m_5=m_2-m_6=m_3=m_4=0.
\end{equation}
On a trouv\'e ainsi les \'equations de $X\times\hild$ dans $\hilt\times\hild$ au voisinage de $(s,b)$. D'apr\`es la remarque \ref{obs17} les \'equations du ferm\'e $\hiltd$ dans $\hilt\times\hild$ sont $[M_{23}]=[M_{13}]=[M_{12}]=0$.
On exprime $[M_{23}],[M_{13}],[M_{12}]$ dans la base donn\'e par $[UW],[W^2]$:
\begin{eqnarray*}
{[M_{23}]}&=&[UW](-m_2-m_6+\alpha-m_3\lambda)+[W^2](\beta-m_3\mu+m_2m_6-m_3m_5)\\
{[M_{13}]}&=&[UW](\lambda\alpha+\mu-m_1+\lambda m_6)+[W^2](\lambda\beta+\mu m_6+m_1m_6-m_3m_4)\\
 {[M_{12}]}&=&[UW](\lambda^2\alpha+2\lambda\mu+m_4+\lambda m_1+\lambda m_5)+[W^2](\lambda^2\beta+\mu^2+\mu m_1+\mu m_5+m_1m_5-m_2m_4).
\end{eqnarray*}

Par cons\'equent les \'equations de $\hiltd$ dans $\hilt\times\hild$ sont ces $6$ coefficients. Les diff\'erentielles des quatre d'entre eux
\begin{eqnarray}
\label{ecuatia29}
-m_2-m_6+\alpha-m_3\lambda\nonumber\\
\lambda\alpha+\mu-m_1+\lambda m_6\nonumber\\
\lambda^2\alpha+2\lambda\mu+m_4+\lambda m_1+\lambda m_5\\
\beta-m_3\mu+m_2m_6-m_3m_5\nonumber
\end{eqnarray}
sont respectivement $-m_2-m_6+\alpha, \mu-m_1, m_4, \beta$, lin\'eairement ind\'ependantes. Puisque $\hiltd$ est de codimension $4$ dans $\hilt\times\hild$, les \'equations (\ref{ecuatia29}) d\'efinissent $\hiltd$ au voisinage de $(s,b)$. 
%%%
%%%{\underline{\bf Question:}} Pourquoi $\uTor^{\comp[\alpha,\beta,\lambda,\mu,m_i]}_q(\comp[\alpha,\beta,\lambda,\mu,m_i]/_{(\ref{ecuatia28})}, \comp[\alpha,\beta,\lambda,\mu,m_i]/_{(\ref{ecuatia29})})$ est nul pour $i>1$ et inversible sur $\comp[\alpha,\beta,\lambda,\mu,m_i]/_{(\ref{ecuatia28})+(\ref{ecuatia29})}$ lorsque $i=1$ ?

On consid\`ere les faisceaux  $F$ et $L$ sur $\hilt\times\hild$  libres dans le voisinage de $(s,b)\in\hilt\times\hild$, de rangs respectifs $3$ et $1$, et on prend les sections locales 
\begin{equation}
\label{ecuatia3001}
f=(m_1-m_5,m_2-m_6,m_3),
\end{equation}
et $l=m_4$ dans $F^*$ respectivement $L^*$. D'apr\`es la pr\'esentation (\ref{ecuatia28}) le faisceau $\maO_{X\times\hild}$ a une r\'esolution de Koszul $\Lambda^{\bullet}(F\oplus L)\to\maO_{X\times\hild}$ sur $\hilt\times\hild$. En tenant compte du fait que $\Lambda^iL=0$ pour $i\ge 2$, cette r\'esolution $\Lambda^{\bullet}(F\oplus L)$ s'\'ecrit:
\begin{eqnarray*}
0\to L\otimes\Lambda^3F\stackrel{ \left(\begin{array}{c}f\\e\end{array}\right) }{\to}L\otimes\Lambda^2F\oplus\Lambda^3F\stackrel{ \left(\begin{array}{cc}f&0\\e&f\end{array}\right) }{\to}L\otimes F\oplus\Lambda^2F&\stackrel{ \left(\begin{array}{cc}f&0\\e&f\end{array}\right) }{\to}\\\
\to L\oplus F\stackrel{ \left(\begin{array}{cc}e&f\end{array}\right) }{\to}\maO_{\hilt\times\hild}\to 0.
\end{eqnarray*}
Le complexe $\Lambda^{\bullet}(F\oplus L)$ rentre dans une suite exacte de complexes:
$$
0\to\Lambda^{\bullet}(F)\to\Lambda^{\bullet}(F\oplus L)\to L\otimes\Lambda^{\bullet}(F)[1]\to 0,
$$
d'o\`u l'exactitude de la suite exacte longue:
\begin{eqnarray}
\cdots\to\has_q(\Lambda^{\bullet}(F)\otimes_{\maO_{\hilt\times\hild}}\maO_{\hiltd})\to\uTor_q^{\maO_{\hilt\times\hild}}(\maO_{\hiltd},\maO_{X\times\hild})\to\nonumber\\ 
\label{ecuatia3002}
 \to L\otimes\has_{q-1}(\Lambda^{\bullet}(F)\otimes_{\maO_{\hilt\times\hild}}\maO_{\hiltd})\to\cdots
\end{eqnarray}
Les \'equations (\ref{ecuatia3001}) et les diff\'erentielles des \'equations (\ref{ecuatia29}) sont lin\'eairement ind\'ependantes. Puisque $Y$ est de dimension $3$ et les \'equations (\ref{ecuatia3001}), (\ref{ecuatia29}) s'annulent sur $Y$, on trouve que $Y$ est d\'efini par ces $7$ \'equations. Donc
$$\Lambda^{\bullet}(F)\otimes_{\maO_{\hilt\times\hild}}\maO_{\hiltd}$$
 est une r\'esolution pour $\maO_Y$ et la suite (\ref{ecuatia3002}) se r\'eduit \`a:
$$0\to\uTor_1^{\maO_{\hilt\times\hild}}(\maO_{\hiltd},\maO_{X\times\hild})\to L\vert_Y\stackrel{l}{\to}\maO_Y\to\uTor_0^{\maO_{\hilt\times\hild}}(\maO_{\hiltd},\maO_{X\times\hild})\to 0.$$
L'\'equation $l=m_4$ est v\'erifi\'ee sur $Y$. 

Alors $\uTor_0^{\maO_{\hilt\times\hild}}(\maO_{\hiltd},\maO_{X\times\hild})$, $\uTor_1^{\maO_{\hilt\times\hild}}(\maO_{\hiltd},\maO_{X\times\hild})$ sont localement libres de rang $1$ au voisinage du point $(s,b)\in Y$. Par l'homog\'en\'eit\'e de l'action du groupe $\pgl(3)$, on obtient l'\'enonc\'e global.

{\rm (iii)} On a d\'emontr\'e que $\uTor_q^{\maO_{\hilt\times\hild}}(\maO_{\hiltd},\maO_{X\times\hild})$ est respectivement $\maO_Y$ pour $q=0$, $M$ pour $q=1$ et $0$ pour $q>0$. On obtient l'\'egalit\'e dans le groupe de K-th\'eorie $\ka(X\times\hild)$:
\begin{equation}
\label{ecuatia500}
[\maO_Y]-[M]=(i\times id)^*[\maO_{\hiltd}].
\end{equation}
o\`u $(i\times id)$ est le \mor $X\times\hild\to\hilt\times\hild$ et $[\maO_{\hiltd}]$ la classe du faisceau $\maO_{\hiltd}$ dans le groupe de K-th\'eorie $\ka(\hilt\times\hild)$.

Puisque le \mor $j:Y\to X\times\hild$ est une immersion ferm\'ee, on peut r\'e-\'ecrire l'\'egalit\'e (\ref{ecuatia500}) sous la forme
\begin{equation}
\label{ecuatia501}
j_{!}(1-[M])=(i\times id)^*[\maO_{\hiltd}].
\end{equation}
On applique le \mor $\ka(X\times\hild)\stackrel{ch}{\to}\has^*(X\times\hild)$ \`a l'\'egalit\'e 
(\ref{ecuatia501}). On utilise la formule de Riemann-Roch pour le membre de gauche et la fonctorialit\'e des classes de Chern pour le membre de droite. On obtient l'\'egalit\'e dans $\has^*(X\times\hild)$:
\begin{equation}
\label{ecuatia502}
j_{*}\left(\frac{1-e^{c_1(M)}}{td N_{Y/X\times\hild}}\right)=(i\times id)^*ch(\maO_{\hiltd}).
\end{equation}
On applique \`a l'\'egalit\'e (\ref{ecuatia502}) la projection $\has^*(X\times\hild)\stackrel{p_*}{\to}\has^{*-4}(X)$ et on consid\`ere l'\'egalit\'e obtenue dans $\has^0(X)$. Le diagramme cart\'esien:
$$\diagram
X\times\hild\rto^{i\times id}\dto_{p}&\hilt\times\hild\dto^{pr}\\
X\rto^{i}&\hilt
\enddiagram$$
implique $p_*(i\times id)^*=i^*pr_*.$ On obtient l'\'egalit\'e dans $\has^0(X)$:
$$-p_*(c_1(M))=i^*(pr_*([\hiltd])).$$
Le \mor $\hiltd\to\hilt$ est fini de degr\'e $3$. Alors $p_*(c_1(M))=-3$, ce qu'il fallait d\'emontrer.\hfill
$\Box$
%\end{itemize}

\vs\par{\bf Preuve du lemme \ref{lema19}:}

L'application $M=(\pi,\tmu):\Sigma\to\hilt\times X$ a son image dans la vari\'et\'e d'incidence $\Xi=X^{\mbox{}^{[3,1]}}\subset\hilt\times X$. On obtient un diagramme commutatif:
$$\diagram
\Sigma\rto^{M}\drto_{\pi}&\Xi\rto^{p_{31}}\dto^{\pi_{31}}&X\\
&\hilt&.
\enddiagram$$
Il en r\'esulte un \mor de faisceaux sur $X$:
$$\eldtr=\pi_{31*}p_{31}^*\eld\to\pi_{31*}(M_*M^*)p_{31}^*\eld=(\pi_{31*}M_*)(M^*p_{31}^*\eld)=\pi_*\tmu^*\eld.$$
Puisque $\tmu$ co\"{\i}ncide avec $\mu:\dhilts\to X$ au-dessus de $\dhilts$, le \mor obtenu prolonge le \mor de la suite exacte (\ref{ecuatia000}). \hfill $\Box$

\vs\par{\bf Preuve du lemme \ref{lema20}:}

On construit le \mor (\ref{ecuatia301,7}) en utilisant le diagramme (\ref{ecuatia302,5}). Le \mor $\tpi_{32}$ v\'erifie par construction $\tpi_{32*}\maO_D=\maO_{\Sigma}.$ Le faisceau $\tmu^*\eld$ est localement libre sur $\Sigma$. On obtient $\tpi_{32*}\tpi^*_{32}\tmu^*\eld=\tmu^*\eld.$
Alors:
$$\pi_*\tmu^*\eld=\pi_*\tpi_{32*}\tpi^*_{32}\tmu^*\eld=\pi_{32*}\ttmu\eld=\pi_{32*}p_{32}^*HC^*\eld.$$
Le faisceau $\eld$ sur $X$ est la restriction de la diagonale $X$ au faisceau $\dedoi$ sur $\esxd$. Alors $HC^*\eld=\dedoi\vert_{\dhild}.$ Par cons\'equent:
\begin{equation}
\label{ecuatia301,9}
\pi_*\tmu^*\eld=\pi_{32*}(p_{32}^*\dedoi\vert_D).
\end{equation}
Par d\'efiniton $\dedoit=\pi_{32*}p_{32}^*\dedoi$ pour les morphismes du diagramme (\ref{ecuatia0}). Le \mor (\ref{ecuatia301,7}) est le \mor canonique
\begin{equation}
\label{ecuatia860}
\pi_{32*}(p_{32}^*\dedoi)\to\pi_{32*}(p_{32}^*\dedoi\vert_D).
\end{equation}
Puisque $\tmu:D\to X$ co\"{\i}ncide avec $\mu:\dhilts\to X$ au-dessus de $\dhilts$, le \mor (\ref{ecuatia301,7}) prolonge le \mor de la suite exacte (\ref{ecuatia000}).
Pour montrer la surjectivit\'e, on rappelle qu'on a d\'emontr\'e  dans la preuve de la proposition \ref{prop14} que le \mor $\pi_{32}$ \'etait fini au-dessus de $\hilt\setminus X$. Par cons\'equent le \mor (\ref{ecuatia860}) est surjectif. \hfill $\Box$

\section{ Calcul de $\has^*(\hilt,\esdoi\eltr)$}
\label{sectiune4}

On calculera dans cette section  $\has^*(\hilt,\esdoi\eltr)$ \`a partir de la suite exacte du \thm \ref{teorema13}. On commence par calculer la cohomologie des termes concern\'es.

\begin{prop}
\label{prop21}

On a
\begin{equation}
\label{ecuatia21,0} 
\has^*(\hilt,\eldtr)=\esdoi\has^*(\maO_X)\otimes\has^*(\eld).
\end{equation}
\end{prop}

\vs\par{\bf Preuve:}

C'est le  r\'esultat (\ref{ecuatia002bis}) pour $k=1$, $n=3$, $A=\maO$. \hfill $\Box$

\begin{prop}
\label{prop22}

On a 
$$\has^*(\hilt,\dedoit)=\has^*(X,\maO_X)\otimes\has^*(\hild,\dedoi).$$
\end{prop}

\vs\par{\bf Preuve:}

On rappelle que $\dedoit=\pi_{32*}\eldoi$, $\eldoi=p_{32}^*\dedoi$, o\`u $\pi_{32}, p_{32}$ sont les morphismes du diagramme (\ref{ecuatia0}). On a d\'emontr\'e que les fibres du \mor $\pi_{32}$ sont finies au-dessus de l'ouvert $\hilt\setminus X$, o\`u $X$ param\`etre les sch\'emas du type $\{m_x^2\}$. On a d\'emontr\'e dans le lemme \ref{lema16}(i) que les fibres du \mor $\pi_{32}$ au-dessus de $X$ sont $\proj_1$ et dans le lemme \ref{lema15} que la restriction du \fib inversible $\eldoi$ \`a ces fibres est le \fib inversible trivial. Alors $\er^q\pi_{32*}\eldoi=0$ pour $q>0$. D'apr\`es la suite spectrale de Leray on obtient:
\begin{equation}
\label{ecuatia404}
\has^q(\hilt,\dedoit)=\has^q(\hiltd,\eldoi).
\end{equation}
On consid\`ere le \mor $\tp_{32}=(a,p_{32}):\hiltd\to X\times\hild$, o\`u le \mor $a$ associe au couple $(Z,Z')\in\hilt\times\hild$ le point $HC(Z)-HC(Z')$, $HC$ \'etant le \mor de Hilbert-Chow. Mais $\eldoi=\tp_{32}^*(\maO\boxtimes\dedoi)$.  On a d\'emontr\'e dans \cite{Danila} que
$$\begin{array}{ccc}
   \er^q\tp_{32*}\maO_{\hiltd}&=&\begin{cases}0 \mbox{ \ \ \ si \ \ \ } q>0\cr \maO_{X\times\hild} \mbox{ \ \ \ si \ \ \ } q=0.\cr\end{cases}
 \end{array}$$
D'apr\`es la suite spectrale de Leray on obtient
\begin{equation}
\label{ecuatia405}
\has^q(\hiltd,\eldoi)=\has^q(X\times \hild,\maO\boxtimes\dedoi).
\end{equation}
Les \'egalit\'es (\ref{ecuatia404}), (\ref{ecuatia405}) et le \th de K\"unneth impliquent le r\'esultat.
 \hfill $\Box$

\begin{prop}
\label{prop23}

On a 
\begin{equation}
\label{ecuatia407}
\has^*(\hilt,\pi_*\tmu^*\eld)=\has^*(X,\maO_X)\otimes\has^*(X,\eld).
\end{equation}
\end{prop}

\vs\par{\bf Preuve:}

Soit $HC:\hilt\to\es^3X$ le \mor de Hilbert-Chow. On consid\`ere le \mor $\nu=(a,\tmu):\Sigma\to X\times X$, o\`u $a$ est le \mor $\Sigma\to X$ donn\'e par $a(s)=HC(\pi(s))-2\tmu(s).$ On a
$$\pi_*\tmu^*\eld=\pi_*\nu^*(\maO\boxtimes\eld).$$
On d\'emontre:
\begin{lema}
\label{lema400}
On a:
$$\begin{array}{ccc}
   \er^q\nu_*\maO_{\Sigma}&=&\begin{cases}0 \mbox{ \ \ \ si \ \ \ } q>0\cr \maO_{X\times X} \mbox{ \ \ \ si \ \ \ } q=0.\cr\end{cases}
 \end{array}$$
\end{lema}

En appliquant le lemme on trouve:
\begin{equation}
\label{ecuatia310}
{\begin{array}{ccc}
   \er^q\nu_*\nu^*(\maO\boxtimes\eld)&=&\begin{cases}0 \mbox{ \ \ \ si \ \ \ } q>0\cr \maO\boxtimes\eld \mbox{ \ \ \ si \ \ \ } q=0.\cr\end{cases}
 \end{array}}
\end{equation}
Alors:
$$\has^*(\hilt,\pi_*\nu^*(\maO\boxtimes\eld))=\has^*(\Sigma,\nu^*(\maO\boxtimes\eld))=\has^*(X\times X,\maO\boxtimes\eld)=\has^*(X,\maO)\otimes\has^*(X,\eld),$$
en appliquant le fait que le \mor $\pi$ est fini, la relation (\ref{ecuatia310}) et le \th de K\"unneth.
\hfill $\Box$

\vs\par{\bf Preuve du lemme \ref{lema400}:}

Consid\'erons le diagramme 
$$\diagram
D\rto^{\tnu}\dto_{\tpi_{32}}&X\times\dhild\dto^{b}\\
\Sigma\rto^{\nu}\dto_{\pi}&X\times X\\
\dhilt&
\enddiagram$$
o\`u $\tnu=(a\circ\tpi_{32},p_{32})$, $b=(id, pr)$, et $pr:\dhild\to X$ est le \mor qui associe au sch\'ema double son support.
On a vu dans la preuve de la proposition \ref{prop14} et dans le lemme \ref{lema16}(i) que les fibres du \mor $\pi_{32}=\pi\circ\tpi_{32}$ sont des ensembles finis ou l'espace projectif $\proj_1$. Par suite les fibres du \mor $\tpi_{32}$ sont des points ou l'espace projectif $\proj_1$. Donc $\er^q\tpi_{32}\maO_D=0$ pour $q>0$. D'apr\`es la construction de $\Sigma$ de la preuve du lemme \ref{lema300} on obtient $\tpi_{32*}\maO_D=\maO_{\Sigma}$. 
La suite spectrale de Leray pour $\er(\nu\circ\tpi_{32})_*=\er\nu_*\circ\er\tpi_{32*}$ implique:
\begin{equation}
\label{ecuatia401}
\er^q\nu_*\maO_{\Sigma}=\er^q(\nu\circ\tpi_{32})_*\maO_D, \mbox{ \ \ \ pour\ \ \ } q\ge 0.
\end{equation}
Le \mor $\tnu$ est un \iso au-dessus des points $(X,Z')$ avec $x\ne\supp Z'$ (son inverse est donn\'e par $(x,Z')\mapsto(\maO_x\oplus\maO_{Z'},Z')$). 
On utilise le

\begin{lema}
\label{lema401}
Soit $Z'\in\dhild$ de support $x$. La fibre $F(x,Z')$ du \mor $\tnu$ dans $(x,Z')$ est une courbe lisse et rationnelle. 
\end{lema}

Alors $\er^q\tnu_*\maO_D=0$ pour $q>0$. Le \mor $\tnu$ est birationnel et la vari\'et\'e $X\times\dhild$ est lisse, donc normale. D'apr\`es Zariski's Main Theorem (\cite{Hart}, III 11.4) on a $\tnu_*\maO_D=\maO_{X\times\dhild}$. La suite spectrale de Leray pour $\er(b\circ\tnu)_*=\er b_*\circ\er\tnu_*$ nous donne
\begin{equation}
\label{ecuatia402}
\er^q(b\circ\tnu)_*\maO_D=\er^q b_*\maO_{X\times\dhild}.
\end{equation}
Finalement $b$ est une fibration \`a fibres $\proj_1$. Alors
\begin{equation}
\label{ecuatia403}
{\begin{array}{ccc}
   \er^qb_*\maO_{X\times\dhild}&=&\begin{cases}0 \mbox{ \ \ \ si \ \ \ } q>0\cr \maO_{X\times X} \mbox{ \ \ \ si \ \ \ } q=0.\cr\end{cases}
 \end{array}}
\end{equation}
Les relations (\ref{ecuatia401}), (\ref{ecuatia402}) et (\ref{ecuatia403}) nous aident \`a conclure la preuve du lemme \ref{lema400}. \hfill $\Box$

\vs\par{\bf Preuve du lemme \ref{lema401}:}

C'est un \'enonc\'e relatif \`a $\maO_X/m_x^3$, donc on peut supposer $X=\pp$. On note $(U:V:W)$ les coordonn\'ees homog\`enes sur $\pp$ et $u=\frac{U}{W}, v=\frac{V}{W}$.  Par homog\'en\'eit\'e on peut supposer $x=(0:0:1)$ et $\maI_{Z'}=(u^2,v)$.
Les sch\'emas de longueur $3$ \`a support dans $x$  contenant $Z'$ sont ceux d'id\'eal $(v-Au^2,u^3), A\in\comp$ ou $(u^2,uv,v^2).$ Il r\'esulte que la fibre g\'eom\'etrique $F(x)$ est ${\mathbb{A}}^1\cup\{*\}.$ Pour d\'emontrer que $F(x)$ est sch\'ematiquement $\proj_1$, il suffit de d\'emontrer que $F(x)$ est lisse au voisinage de chacun de ses points.
Les points d'id\'eal $(v-Au^2,u^3)$, $A\ne 0$, sont dans la m\^eme orbite de l'action de $\pgl(3)$. Tout voisinage du sch\'ema $(v,u^3)$ contient des points $(v-Au^2,u^3)$. Il suffit donc de d\'emontrer la lissit\'e dans le voisinage des points d'id\'eal $(v,u^3)$ et $(u^2,uv,v^2)$.

On utilise la param\'etrisation du lemme \ref{lema3.600}. Les sch\'emas du voisinage du point $(v,u^3)$ de $\hilt$ qui contiennent le sch\'ema $(v,u^2)$, v\'erifient les \'equations (\ref{ecuatia3078}) avec $\lambda=\mu=\alpha=0.$ Ils sont donc les points d'\'equation $B=C=q=r=0$, c'est-\`a-dire les sch\'emas $(v-Au^2,u^3+pu^2)$. L'image par $\tnu$ d'un tel sch\'ema dans $X\times\dhild$ est $((-p,-Ap^2),Z')$. Pour que l'image soit $(x,Z')$, on doit avoir $p=Ap^2=0$. Alors $F(x,Z')$ est d'\'equation $B=C=p=q=r=0$ au voisinage de du point $(v,u^3)$, donc elle est lisse.

Dans le voisinage du point $(U^2,UV,V^2)=(u^2,uv,v^2)$ on utilise la param\'etrisation du lemme \ref{lema3.11}. Les sch\'emas du voisinage qui contiennent le sch\'ema $(V,U^2)=(v,u^2)$ v\'erifient les \'equations (\ref{ecuatia29}) avec $\lambda=\mu=\alpha=\beta=0$, c'est-\`a-dire:
$$m_2+m_6=m_1=m_4=m_2m_6-m_3m_5=0.$$
Elles sont en correspondance avec les sch\'emas d'id\'eal $(v(v+m_5),v(u+m_2),u^2-vm_3).$ L'image par $\tnu$ d'un tel sch\'ema est $((-m_2,-m_5),Z')$. Alors 
$F(x,Z')$ est d'\'equations $m_1=m_2=m_4=m_5=m_6$ au voisinage du point 
$(u^2,uv,v^2)$, donc elle est lisse. 
 \hfill $\Box$

%%\begin{teor}
%%\label{teor24}
%%Soit $X$ une surface projective lisse et $L$ un faisceau inversible sur $X$. L'application canonique (\ref{ecuatia18})  pour $n=3$:
%%\begin{equation}
%%\label{ecuatia31}
%%can:\esdoi\has^*(X,L)\to\has^*(\hilt,\esdoi\eltr),
%%\end{equation}
%%induit la d\'ecomposition en somme directe
%%$$\has^*(\hilt,\esdoi\eltr)=\has^*(\maO_X)\otimes\esdoi\has(X,L)\bigoplus(\esdoi \has^*(\maO_X)/{\has^*(\maO_X)})\otimes\has(X,\eld).$$
%%%
%%%
%%%Elle est injective et a pour conoyau:
%%%$$(\esdoi \has^*(\maO_X)/{\has^*(\maO_X)})\otimes\has(X,\eld)\bigoplus(\has^*(\maO_X)/{\comp})\otimes\esdoi\has^*(L).$$
%%\end{teor}

\vs On rappelle la convention faite dans la section \ref{sectiune2} d'omettre  l'espace $X$ dans la notation $ \has^*(X,L)$ pour tout faisceau inversible $L$.

\vs\par{\bf Preuve du \th \ref{teor24}:}

Le \th \ref{teorema13} nous conduit \`a une suite longue de cohomologie:
$$\cdots\to\has^*(\hilt,\esdoi\eltr)\to\has^*(\hilt,\dedoit)\oplus\has^*(\hilt,\eldtr)\to\has^*(\hilt,\pi_*\tmu^*\eld)\to\cdots$$
D'apr\`es les propositions \ref{prop21}, \ref{prop22}, le lemme \ref{lema9} et la proposition  \ref{prop23} on peut r\'e-\'ecrire cette suite sous la forme:
\begin{equation}
\label{ecuatia405,5}
\cdots\to\has^*(\hilt,\esdoi\eltr)\to\has^*(\maO)\otimes\esdoi\has^*(L)\bigoplus\esdoi\has^*(\maO)\otimes\has^*(\eld)\to\has^*(\maO)\otimes\has^*(\eld)\to\cdots
\end{equation}

\begin{lema}
\label{lema}
Le diagramme
$$\diagram
\has^*(\maO)\otimes\esdoi\has^*(L)\dto_{can}\drto^{id}&\\
\has^*(\hilt,\esdoi\eltr)\rto_{u\ \ \ \ \ \ \ }&\has^*(\dedoit)=\has^*(\maO)\otimes\esdoi\has^*(L),
\enddiagram$$
%%%o\`u $w$ est donn\'e par $\alpha\to 1\otimes\alpha$, 
est commutatif.
\end{lema}

\vs\par{\bf Preuve:}

Le \mor $\esdoi\eltr\to\dedoit$ s'\'ecrit $\pi_*\maO_{\esdoi(\Xi)}\to\pi_{32*}\maO_{\hiltd}$ pour $L=\maO$, dans les notations du diagramme (\ref{ecuatia44,5}). On obtient un \mor $U:\hiltd\to\esdoi(\Xi)$ qui induit en cohomologie l'application $u$.

L'application $can$ est induite par d\'efinition par le \mor $P:\esdoi(\Xi)\to X\times\esxd$.

L'application $id$ est induite, d'apr\`es la proposition \ref{prop22}, par le \mor $\hiltd\stackrel{(a,HC\circ p_{32})}{\longrightarrow}X\times\esxd.$
Il suffit de d\'emontrer que le diagramme
$$\diagram
X\times \esxd&\\
\esdoi(\Xi)\uto^{P}&\hiltd\ulto_{(a,HC\circ p_{32})}\lto_{U}
\enddiagram$$
est commutatif. Il suffit de le v\'erifier au-dessus de l'ouvert dense $\pi_{32}^{-1}(\hiltss)$, o\`u $\hiltss$ est l'ouvert des sch\'emas lisses de $\hilt$. Ici l'application $U$ est d\'efinie par:
$$((x_1,x_2,x_3),(x_i,x_j))\mapsto((x_1,x_2,x_3),(x_i+x_j)),\ \ i\ne j,\ i,j\in\{1,2,3\}.$$
La commutativit\'e du diagramme en est une cons\'equence.
 \hfill $\Box$

\vs Le lemme implique la d\'ecomposition en somme directe
$$\has^*(\hilt,\esdoi\eltr)=\has^*(\maO)\otimes\esdoi\has^*(L)\oplus K^*,$$
et la suite (\ref{ecuatia405,5}) peut s'\'ecrire:
$$\cdots K^*\to \esdoi\has^*(\maO)\otimes\has^*(\eld)\to \has^*(\maO)\otimes\esdoi\has^*(\eld)\to K^{*+1}\to\cdots$$

Afin de comprendre le \mor $\esdoi\has^*(\maO)\otimes\has^*(\eld)\stackrel{b}{\to}\has^*(\maO)\otimes\has^*(\eld)$ on consid\`ere le diagramme:
\begin{equation}
\label{ecuatia406}
{\diagram
\Xi\rto^{\ \ c\ \ \ \ \ }&\esxd\times X\\
\Sigma\uto^{(\pi,\tmu)}\rto_{\ \ \nu=(a,\tmu)\ \ \ \ \ }&X\times X\uto_{d},
\enddiagram}
\end{equation}
o\`u:

$\bullet$  $\Xi=X^{\mbox{}^{[3,1]}}$ est la vari\'et\'e d'incidence d\'efinie \`a l'aide du diagramme (\ref{ecuatia00})

$\bullet$ le \mor $c$ associe au couple $(Z,x)\in\hilt\times X$ le point $(HC(Z)-x,x)\in\esxd\times X$

$\bullet$ le \mor $d$ est d\'efini par $d(x,y)=(x+y,y).$

Le diagramme (\ref{ecuatia406}) est commutatif par d\'efinition  des morphismes. On a d\'emontr\'e dans \cite{Danila} que l'\iso (\ref{ecuatia21,0}) provient de l'identification:
$$\has^*(\hilt,\eldtr)=\has^*(\Xi,c^*(\maO\boxtimes\eld))=\has^*(\esxd\times X, \maO\boxtimes\eld).$$
On a d\'emontr\'e dans la proposition \ref{prop23} que l'\iso (\ref{ecuatia407}) provient de l'identification:
\begin{equation}
\label{ecuatia408}
\has^*(\hilt,\pi_*\tmu^*\eld)=\has^*(\Sigma,\nu^*(\maO\boxtimes\eld))=\has^*(X\times X, \maO\boxtimes\eld).
\end{equation}
En utilisant le diagramme commutatif 
(\ref{ecuatia406}) et les identifications (\ref{ecuatia407}), (\ref{ecuatia408}) on d\'eduit que le \mor $\esdoi\has^*(\maO)\otimes\has^*(\eld)\to\has^*(\maO)\otimes\has^*(\eld)$ est le \mor
$$d^*:\has^*(\esxd\times X,\maO\boxtimes\eld)\to\has^*(X\times X, \maO\boxtimes\eld).$$
D'apr\`es la formule de K\"unneth il s'\'ecrit:
$$uv\otimes\alpha\to u\otimes v\alpha+(-1)^{pq}v\otimes u\alpha,$$
pour $u\in\has^p(\maO_X)$, $v\in\has^q(\maO_X)$, $\alpha\in\has^*(\eld).$
En particulier ce \mor est surjectif: le \mor 
\begin{equation}
\label{ecuatia408,5}
\has^*(\maO)\otimes\has^*(L)\to\esdoi\has^*(\maO)\otimes\has^*(L)
\end{equation}
donn\'e par
$$u\otimes\alpha\mapsto (1\cdot u)\otimes\alpha-\frac{1}{2}1\otimes u\alpha$$
est une section de ce morphisme.
Alors $K^*=\Ker b$, d'o\`u la conclusion.
%%%
%%%
%%%Par cons\'equent la suite exacte longue (\ref{ecuatia405,5}) se s\'epare en suites exactes courtes:
%%%\begin{equation}
%%%\label{ecuatia409}
%%%0\to\has^*(\hilt,\esdoi\eltr)\to\has^*(\maO)\otimes\esdoi\has^*(L)\bigoplus\esdoi\has^*(\maO)\otimes\has^*(\eld)\to\has^*(\maO)\otimes\has^*(\eld)\to 0.
%%%\end{equation}
%%%
%%%
%%%\vs D'apr\`es ce lemme le \mor $can$ est injectif et on peut r\'e-\'ecrire la suite exacte (\ref{ecuatia409}) sous la forme:
%%%\begin{eqnarray}
%%%\label{ecuatia410}
%%%0\to\esdoi\has^*(L)\to\has^*(\hilt,\esdoi\eltr)\to(\has^*(\maO)\otimes\esdoi\has^*(L))/\esdoi\has^*(L)\bigoplus\esdoi\has^*(\maO)\otimes\has^*(\eld)\to\hfill\\
%%%\hfill\to\has^*(\maO)\otimes\has^*(\eld)\to 0.\nonumber
%%%\end{eqnarray}
%%%La d\'ecomposition en somme directe
%%%$$\has^*(\maO)\otimes\esdoi\has^*(L)=\esdoi\has^*(L)\bigoplus(\has^*(\maO)/\comp)\otimes\esdoi\has^*(L)$$
%%%et la d\'ecomposition en somme directe donn\'ee par la section (\ref{ecuatia408,5}) nous aident \`a r\'e-\'ecrire la suite sous la forme \'enonc\'ee dans le th\'eor\`eme. 
\hfill $\Box$

%%\begin{cor}
%%\label{cor25}
%%Soit $X$ une surface projective lisse et $L$ un faisceau inversible sur $X$. Si $X$ satisfait en outre $q=p_g=0$, le morphisme canonique (\ref{ecuatia31}) est un isomorphisme.
%%
%%\end{cor}

\vs\par{\bf Preuve du \th \ref{teor4.6'}:}

On tensorise la suite exacte du \th \ref{teorema13} par le \fs \inve $\dta$ sur $\hilt$. De la m\^eme mani\`ere que les \'enonc\'es \ref{prop21}, \ref{prop22}, \ref{prop23} on prouve
$$\has^*(\hilt,\eldtr\otimes\dta)=\esdoi\has^*(X,A)\otimes\has^*(X,\eld\otimes A),$$
$$\has^*(\hilt,\dedoit\otimes\dta)=\has^*(X,A)\otimes\esdoi\has^*(X,L\otimes A),$$ 
$$\has^*(\hilt,\pi_*\tmu^*\eld\otimes\dta)=\has^*(X,A)\otimes\has^*(X,\eld\otimes A^2).$$
On obtient la suite exacte longue:
\begin{eqnarray}
\cdots\to \has^*(\hilt,\esdoi\eltr\otimes\dta)\to \has^*(X,A)\otimes\esdoi\has^*(X,L\otimes A)\oplus \esdoi\has^*(X,A)\otimes\has^*(X,\eld\otimes A)\to\nonumber\\
\label{ecuatia3002b}
\to\has^*(X,A)\otimes\has^*(X,\eld\otimes A^2)\to\cdots
\end{eqnarray}
De m\^eme que dans la preuve du lemme \ref{lema}, la composition
$$\has^*(X,A)\otimes\esdoi\has^*(X,L\otimes A)\stackrel{can}{\to}\has^*(\hilt,\esdoi\eltr\otimes\dta)\to \has^*(X,A)\otimes\esdoi\has^*(X,L\otimes A)$$
est un isomorphisme. Alors $\has^*(\hilt,\esdoi\eltr\otimes\dta)=\has^*(X,A)\otimes\esdoi\has^*(X,L\otimes A)\oplus K^*$, et la suite exacte (\ref{ecuatia3002b}) se transforme dans la suite exacte (\ref{ecuatia3001b}). \hfill $\Box$

\section{ Le calcul de $\has^0(\hilx,\esdoi\elen)$}
\label{sectiune5}

Cette section est consacr\'ee au calcul de l'espace des sections globales  $\has^0(\hilx,\esdoi\elen)$. 
On suppose partout dans cette section que $n\ge 2$.
Puisque $\hilxs$ est un ouvert dont la codimension du compl\'ementaire est \'egale \`a $2$ dans la vari\'et\'e lisse $\hilx$, et $\esdoi\elen$ est localement libre, on a:
\begin{equation}
\label{ecuatia2501}
\has^0(\hilx,\esdoi\elen)=\has^0(\hilxs,\esdoi\elen).
\end{equation}
On calculera $\has^0(\hilxs,\esdoi\elen)$ \`a partir de la suite exacte (\ref{ecuatia000}). On commence par calculer $\has^0(\hilxs,\dedoin)$ et $\has^0(\hilxs,\mu^*\eld\vert_{\dhilxs})$. 

On note $L_i=pr_i^*L$ l'image r\'eciproque de $L$ par la projection $pr_i:\xn\to X$ et $L_{ij}=L_i\otimes L_j$. Par abus de notation on note \'egalement $L_{ij}$ l'image r\'eciproque $\rho^*L_{ij}$ sur $\bens$, $\rho$ \'etant le \mor du diagramme (\ref{ecuatia900}).

\begin{prop}
\label{prop151}
Le faisceau $\dedoin$ rentre dans la suite exacte sur $\bens$:
\begin{equation}
\label{ecuatia2502}
0\to\dedoin\to\sum_{i\ne k}L_{ik}\stackrel{a}{\to}\sum_{\stackrel{i<j}{k\ne i,j}}L_{ik}\vert_{E_{ij}}\to 0,
\end{equation}
o\`u les indices $i,j,k$ parcourent l'ensemble $\{1,2,\ldots,n\}$ et l'application $a$ est d\'efinie de la mani\`ere suivante

$\bullet$ la composante $L_{ik}\to L_{ik}\vert_{E_{ij}}$ est la restriction \`a $E_{ij}$ du faisceau $L_{ik}$ lorsque $i<j$

$\bullet$ la composante $L_{ik}\to L_{ik}\vert_{E_{ij}}=L_{jk}\vert_{E_{ij}}$ est $(-1)$ fois la restriction \`a $E_{ij}$ du faisceau $L_{jk}$ lorsque $i<j$

$\bullet$ la composante $L_{ik}\to L_{ik}\vert_{E_{ij}}$ est nulle dans tous les autres cas.

\end{prop}

\vs\par {\bf Preuve:}

On peut recouvrir $\bens$ par les ouverts $\benij$. Il suffit donc de prouver l'exactitude de la suite exacte en restriction \`a chacun de ces ouverts. On utilise l'identification du lemme \ref{lema3}. Pour simplifier on prend $\{i,j\}=\{1,2\}$. D'apr\`es la relation (\ref{ecuatia7}) on a:
$$\dedoin=\dedoi\boxtimes\maO\oplus\eled\boxtimes(\sum_{i\ge 3}L_i)\oplus\maO\boxtimes(\sum_{3\le i<j\le n}L_{ij}).$$
Par d\'efinition on a:
$$\sum_{i\ne k}L_{ik}=\dedoi\boxtimes\maO\oplus (L_1\oplus L_2)\boxtimes(\sum_{i\ge 3}L_i)\oplus\maO\boxtimes(\sum_{3\le i<j\le n}L_{ij}),$$
$$\sum_{k\ne 1,2}L_{1k}\vert_{E_{12}}=0\oplus L_1\vert_E\boxtimes(\sum_{i\ge 3}L_i)\oplus 0,$$
$$L_{ik}\vert_{E_{ij}}=0 \ \ \mbox{\ \ pour\ \ } \{i,j\}\ne\{1,2\}.$$
La suite exacte (\ref{ecuatia2502}) r\'esulte de la suite exacte (\ref{ecuatia0,7}) sur $\bedoi$:
$$0\to\eled\to L_1\oplus L_2\to L_1\vert_E\to 0.\ \ \hfill \Box$$
%d\'emontr\'ee dans \cite{Danila}.

\begin{prop}
\label{prop152}
Soit $n\ge 2$. On a
$$\has^0(\hilxs,\dedoin)=\esdoi\has^0(X,L).$$
\end{prop}

\vs\par {\bf Preuve:}

Le cas $n=2$ a \'et\'e examin\'e dans le lemme \ref{lema9}. On supposera dans la suite que $n>2$. Le faisceau $\dedoin$ sur $\bens$ est l'image r\'eciproque du faisceau $\dedoin$ sur $\hilxs$. On trouve
$$\has^0(\hilxs,\dedoin)=\has^0(\bens,\dedoin)^{\sigmn}.$$
Le \mor $\rho:\bens\to\xns$ est l'\'eclatement des diagonales $\Delta_{ij}=\{x_i=x_j\}$ de l'ouvert $\xns$ de $\xn$ obtenu en gardant seulement les $n$-uples admettant au plus deux coordonn\'ees \'egales. Alors $\rho_*\maO_{\bens}=\maO_{\xns}$, donc $\has^0(\bens,L_{ik})=\has^0(\xns,L_{ik}).$ Comme $\xns$ est un grand ouvert de $\xn$ (c'est-\`a-dire un ouvert dont le compl\'ementaire est au moins de codimension $2$) on trouve  
$\has^0(\xns,L_{ik})=\has^0(\xn,L_{ik})$. Finalement, la formule de K\"unneth identifie $ \has^0(\xn,L_{ik})=\has^0(X,L)\otimes\has^0(X,L)$, par le \mor
$$\xn\stackrel{(pr_i,pr_k)}{\longrightarrow}X\times X.$$
En conclusion, un \'el\'ement de $\has^0(\bens,L_{ik})$ correspond \`a un \'el\'ement $a_{ik}\in\has^0(X,L)\otimes\has^0(X,L)$. Par construction, l'\'el\'ement $a_{ki}$ est le transpos\'e de l'\'el\'ement $a_{ik}$ par l'action de $\sigmd$ sur $\has^0(X,L)\otimes\has^0(X,L)$.

De mani\`ere analogue l'application $E_{ij}\stackrel{\rho}{\to}\Delta_{ij*}\simeq X^{n-1}_*$ induit un \iso
$$
\has^0(\bens,L_{ik}\vert_{E_{ij}})\simeq \has^0(X,L)\otimes\has^0(X,L),$$
et le \mor $L_{ik}\to L_{ik}\vert_{E_{ij}}$ induit un \iso
$$ \has^0(X,L)\otimes\has^0(X,L)\simeq \has^0(X,L)\otimes\has^0(X,L).$$
\`A partir de la suite exacte (\ref{ecuatia2502}) on trouve qu'un \'el\'ement de $\has^0(\bens,\dedoin)$ s'identifie avec une suite d'\'el\'ements $(a_{ik})_{i\ne k}\in \has^0(X,L)\otimes\has^0(X,L)$ dont les images par l'application $a$ sont nulles.
Pour tout triplet de nombres distincts $i,j,k$  on a:
$$a_{ik}=a_{jk}, a_{ij}=a_{kj}, a_{ji}=a_{ki}.$$
Donc $a_{ik}=a_{jk}=a_{ji}=a_{kj}=a_{ij}$, tous les \'el\'ements $a_{ik}$ sont \'egaux et sym\'etriques: $a_{ik}\in\esdoi\has^0(X,L)\subset  \has^0(X,L)\otimes\has^0(X,L)$.
Alors $\has^0(\bens,\dedoin)=\esdoi\has^0(X,L)$. L'action sur $ \esdoi\has^0(X,L)$ par cette identification est triviale ($\sigmn$ envoie $a_{ik}$ sur $a_{\sigma(i)\sigma(k)}$). On obtient $\has^0(\hilxs,\dedoin)=\esdoi\has^0(X,L)$. \hfill $\Box$

\begin{prop}
\label{prop153}
Soit $n\ge 2$. On a
$$\has^0(\hilxs,\mu^*\eld\vert_{\dhilxs})=\has^0(X,\eld).$$
\end{prop}

\vs\par {\bf Preuve:}

La vari\'et\'e $\dhilxs$ co\"{\i}ncide avec l'ouvert $U$ des points $(Z',x_3+\cdots+x_n)$ dans $\dhild\times\es^{n-2}X$ pour lequels les points $x_i$ sont distincts entre eux et disjoints de $\supp Z'$. On consid\`ere l'ouvert $V\subset X\times\es^{n-2}X$ des points $(x_1,x_3+\cdots+x_n)$ tels que tous les $x_i$ sont distincts.
Le \mor 
$$\pi:U\to V, \ \ (Z',x_3+\cdots+x_n)\mapsto(\supp Z',x_3+\cdots+x_n)$$
est une fibration \`a fibres $\proj_1$, et $\mu=\pi\circ pr_1$. Alors 
$$ \has^0(\hilxs,\mu^*\eld)=\has^0(V,\eld\boxtimes\maO).$$
L'ouvert $V$ est grand dans $X\times\es^{n-2}X$ dans le sens  d\'ej\`a employ\'e, donc
$$  \has^0(V,\eld\boxtimes\maO)=\has^0(X\times\es^{n-2}X,\eld\boxtimes\maO)=\has^0(X,L).\ \hfill \Box$$

%\begin{teor}
%\label{teorema154}
%Le \mor canonique (\ref{ecuatia18}) est un \iso en degr\'e $0$:
%$$
%can:\esdoi\has^0(X,L)\stackrel{\sim}{\to}\has^0(\hilx,\esdoi\elen).
%$$
%\end{teor}

\vs\par {\bf Preuve du \th \ref{teorema154}:}

La suite exacte (\ref{ecuatia000}) induit la suite exacte
\begin{equation}
\label{ecuatia2503}
0\to\has^0(\hilxs,\esdoi\elen)\to\has^0(\hilxs,\dedoin)\oplus\has^0(\hilxs,\elden)\to \has^0(\hilxs,\mu^*\eld\vert_{\dhilxs})
\end{equation}
Le m\^eme argument que dans la preuve du lemme \ref{lema} prouve que la composition
$$\esdoi\has^0(X,L)\stackrel{can}{\to}\has^0(\hilxs,\esdoi\elen)\to\has^0(\hilxs,\dedoin)=\esdoi\has^0(X,L)$$
est un isomorphisme. Par cons\'equent $\has^0(\hilxs,\esdoi\elen)$ s'\'ecrit $\esdoi\has^0(X,L)\oplus K^0$, o\`u $K^0$ est le noyau du morphisme:
$$\has^0(\hilxs,\elden)\stackrel{b}{\to}\has^0(\hilxs,\mu^*\eld\vert_{\dhilxs})=\has^0(X,\eld).$$
Dans \cite{Danila} il est d\'emontr\'e que l'\iso (\ref{ecuatia002bis})
$$\has^0(\hilxs,\elden)=\has^0(X,\eld)$$
est induit par une application $\Xi=X^{\mbox{}^{[n,1]}}\to\es^{n-1}X\times X$. La proposition \ref{prop153} nous donne un \iso
$$\has^0(\dhilxs,\mu^*\eld)=\has^0(X,\eld),$$
induit par un \mor $\dhilx\to\es^{n-2}X\times X$. Par le m\^eme argument que dans la preuve du \th \ref{teor24} (diagramme (\ref{ecuatia406})) ces morphismes sont compatibles. Par cons\'equent le \mor $b$ est un isomorphisme, donc $K^0=0$. 
\hfill $\Box$

\vs\par{\bf Preuve du \th \ref{teor154'}:}

On tensorise la suite exacte (\ref{ecuatia000}) par le faisceau \inve $\dna$ sur $\hilxs$. On obtient une suite exacte de sections globales:
\begin{eqnarray}
0\to\has^0(\hilxs,\esdoi\elen\otimes\dna)&\to&\has^0(\hilxs,\dedoin\otimes\dna)\oplus\has^0(\hilxs,\elden\otimes\dna)\to\nonumber\\
\label{ecuatia2520}
&\to&\has^0(\hilxs,\mu^*\eld\vert_{\dhilxs}\otimes\dna)
\end{eqnarray}
De mani\`ere analogue \`a la preuve des propositions \ref{prop152}, \ref{prop153} on d\'emontre:
\begin{eqnarray}
\has^0(\hilxs,\dedoin\otimes\dna)&=&\es^{n-2}\has^0(A)\otimes\esdoi\has^0(L\otimes A)\nonumber\\
\label{ecuatia2521}
\has^0(\hilxs,\mu^*\eld\vert_{\dhilxs}\otimes\dna)&=&\es^{n-2}\has^0(A)\otimes\has^0(\eld\otimes A^2).
\end{eqnarray}
Le m\^eme argument que dans la preuve du lemme \ref{lema} d\'emontre que la composition 
$$\es^{n-2}\has^0(X,A)\otimes\esdoi\has^0(L\otimes A)\stackrel{can}{\to}\has^0(\hilx,\esdoi\elen\otimes\dna)\to\has^0(\hilxs,\dedoin\otimes\dna)$$
est un isomorphisme. D'o\`u la d\'ecomposition en somme directe (\ref{ecuatia2510}), o\`u $K^0$ est le noyau du \mor
\begin{equation}
\label{ecuatia2522}
b:\has^0(\hilxs,\elden\otimes\dna)\to\has^0(\hilxs,\mu^*\eld\vert_{\dhilxs}\otimes\dna).
\end{equation}
Dans \cite{Danila} il est d\'emontr\'e que l'\iso  (\ref{ecuatia002})
$$\has^0(\hilxs,\elden\otimes\dna)\simeq\has^0(\es^{n-1}X\times X,\maD_{n-1}^A\boxtimes(\eld\otimes A))\simeq\es^{n-1}\has^0(A)\otimes\has^0(\eld\otimes A)$$
est induit par une application $\Xi=\hilxnu\to\es^{n-1}X\times X.$ L'\iso (\ref{ecuatia2521}) est induit par une application $\dhilxs\to \es^{n-2}X\times X.$ Par suite le \mor (\ref{ecuatia2522}) s'\'ecrit sous la forme (\ref{ecuatia2515}). Le m\^eme argument que dans la preuve du \th \ref{teor24} (diagramme (\ref{ecuatia406})) d\'emontre que le \mor (\ref{ecuatia2515}) est induit par l'application:
\begin{eqnarray*}
d:\es^{n-2}X\times X&\to& \es^{n-1}X\times X\\
(x_3+\cdots+x_n,x)&\mapsto&(x+x_3+\cdots+x_n,x).
\end{eqnarray*}
D'apr\`es la formule de K\"unneth l'application (\ref{ecuatia2515}) s'\'ecrit explicitement sous la forme (\ref{ecuatia2516}). \hfill $\Box$

\section{ R\'esultats pour $\has^*(\hilx,\esdoi\elen)$ pour $n$ g\'en\'eral}
\label{sectiune6}

Partout dans cette section $n$ sera un entier $\ge 2$. 
Les r\'esultats de la section \ref{sectiune3} seront utilis\'es ici pour \'etendre la suite exacte (\ref{ecuatia000}) \`a une suite exacte  sur un ouvert $\tihilx$ dont le compl\'ementaire est de codimension $3$ dans  $\hilx$. L'extension \`a $\hilx$ de cette suite exacte, donn\'ee dans le \th \ref{teorema66},  nous permettra d'\'enoncer, dans la proposition \ref{prop67}, des conditions suffisantes pour que le calcul de $\has^*(\hilx,\esdoi\elen)$ pour $n$ g\'en\'eral puisse se faire de mani\`ere analogue que dans les cas particuliers $n=2,3$. 

On note $\trhilx$ l'ouvert des sch\'emas $Z$ dont le support  sch\'ematique est de la forme $HC(Z)=x_1+\cdots+x_n$ ou $2x_1+x_3+\cdots+x_n$ ou $3x_1+x_4+\cdots+x_n$, pour des $x_i$ distincts. On note $\ddhilx$ l'ouvert des sch\'emas $Z$ dont le support sch\'ematique est de la forme $HC(Z)=x_1+\cdots+x_n$ ou $2x_1+x_3+\cdots+x_n$ ou $2x_1+2x_3+x_5+\cdots+x_n$, pour des $x_i$ distincts. On note $\tihilx=\trhilx\cup\ddhilx$.

\begin{prop}
\label{prop61}
Le faisceau $\dedoin$ est localement libre de rang $\frac{n(n-1)}{2}$ sur $\tihilx$.
\end{prop}

\vs\par {\bf Preuve:}

Il suffit de d\'emontrer cette affirmation pour l'image r\'eciproque de $\dedoin$ par le \mor $q:\ben\to\hilx$  du diagramme (\ref{ecuatia900}).
Par abus de notation, on note aussi $\dedoin$ cette image r\'eciproque. 
On note $\tiben, \trben$ respectivement $\ddben$ les images r\'eciproques par $q$ des ouverts $\tihilx, \trhilx, \ddhilx$. On d\'emontrera que $\dedoin$ est localement libre de rang $\frac{n(n-1)}{2}$ sur $\trben$ et $\ddben$. 
L'ouvert $\trben$ est recouvert par les ouverts $\benijk$ des points $Z\in\ben$ tels que $\rho(Z)$ a tous ses points $x_l, l\ne i,j,k$ distincts entre eux et distincts de $x_i, x_j, x_k$. De la m\^eme mani\`ere que dans le lemme \ref{lema3} on peut d\'emontrer que chacun des ouverts $\benijk$ s'identifie \`a un ouvert dans $\be^3\times X^{n-3}$. De la m\^eme mani\`ere que dans le lemme \ref{lema5}, relation (\ref{ecuatia7}), on peut d\'emontrer l'identification (pour simplifier on prend $(i,j,k)=(1,2,3)$):
$$\dedoin=\dedoit\boxtimes\maO\oplus\elded\boxtimes(\sum_{i\ge 4}L_i)\oplus\maO\boxtimes(\sum_{4\le i<j\le n}L_{ij}).$$
D'apr\`es la proposition \ref{prop14}, $\dedoit$ est localement libre de rang $3$. Les autres faisceaux intervenant dans la d\'ecomposition sont localement libres. Par suite $\dedoin$ est localement libre sur $\be^n_{1,2,3}$ de rang $3+3\cdot (n-3)+\frac{(n-3)(n-4)}{2}=\frac{n(n-1)}{2}$.

De la m\^eme mani\`ere, $\ddben$ est recouvert par les ouverts $\benijkl$ pour $i,j,k,l$ distincts dans l'ensemble $\{1,\cdots,n\}$, o\`u $\benijkl$ est l'ouvert des points $Z\in\ben$ tels que $\rho(Z)$ a tous ses points $x_m$ distincts, sauf peut-\^etre $x_i=x_j$ et $x_k=x_l$. Pareil, on prouve que chacun des $\benijkl$ s'identifie \`a un ouvert dans $\bedoi\times\bedoi\times X^{n-4}$
et que dans cette identification on a
(pour simplifier on prend $(i,j,k,l)=(1,2,3,4)$):
\begin{eqnarray*}
\dedoin&=&\dedoi\boxtimes\maO\boxtimes\maO\oplus\maO\boxtimes\dedoi\boxtimes\maO\oplus\elded\boxtimes\elded\boxtimes\maO\oplus\\
&&\oplus\elded\boxtimes\maO\boxtimes (\sum_{i\ge 5}L_i)\oplus\maO\boxtimes\elded\boxtimes(\sum_{i\ge 5}L_i)\oplus \maO\boxtimes\maO\boxtimes(\sum_{5\le i<j\le n}L_{ij}).
\end{eqnarray*}
Par cons\'equent $\dedoin$ est localement libre de rang
 $$1+1+2\cdot 2+2\cdot (n-4)+2\cdot(n-4)+\frac{(n-4)(n-5)}{2}=\frac{n(n-2)}{2}$$
sur $\ddben$. \hfill $\Box$

On avait not\'e $\dhilx$ l'hypersurface des sch\'emas singuliers dans $\hilx$. L'ouvert $\dhilxs=\dhilx\cap\hilxs $ est lisse, irr\'eductible et dont le compl\'ementaire est de codimension $2$ dans $\hilx$. Alors $\dhilx$ est r\'eduite et irr\'eductible, donc int\`egre. On consid\`ere sa normalisation $\pi:\Sigma\to\dhilx$. On note $\tisigma$ l'image r\'eciproque $\pi^{-1}(\dhilx\cap\tihilx)$. L'ouvert $\dhilxs$ est lisse, donc $\pi$ est un \iso au-dessus de $\dhilxs$. On note $i:\dhilxs\to\Sigma$ le \mor d'inclusion. L'analogue du lemme \ref{lema300} est:
\begin{lema}
\label{lema62}
On consid\`ere le diagramme:
\begin{equation*}
{\diagram
\tisigma\dto_{\pi}&&\\
\dhilx&\dhilxs\lto_{j}\ulto_{i}\rto^{\mu}&X,
\enddiagram}
\end{equation*}
o\`u $\mu$ est le \mor d\'efini dans la section \ref{sectiune1}. Il existe une application r\'eguli\`ere $\tmu:\tisigma\to X$ qui rend commutatif ce diagramme.
\end{lema}

\vs\par{\bf Preuve:}

Consid\'erons l'image r\'eciproque 
\begin{equation}
\label{ecuatia600}
\den=p_{n2}^{-1}(\dhild)
\end{equation}
de l'hypersurface $\dhild\subset\hild$. De la m\^eme mani\`ere que dans la d\'emonstration du lemme \ref{lema300}, la preuve du lemme se r\'eduit \`a l'analogue de lemme \ref{lema}:
\begin{lema}
\label{lema63}
La vari\'et\'e $\tiden=\den\cap\pi_{32}^{-1}(\tihilx)$ est lisse.
  \end{lema}

\vs\par{\bf Preuve:}

La question est locale et il suffit de la traiter en g\'eom\'etrie analytique. Le \mor $\pi_{n2}:\den\to\dhilx$ est un \iso au-dessus de $\dhilxs$, et $\dhilxs$ est lisse. Il suffit donc de traiter le probl\`eme au-dessus d'un point de $\trhilx\setminus\hilxs$ et  au-dessus d'un point de $\ddhilx\setminus\hilxs$. On note $Z_0$ ce point.

Dans le premier cas le point $Z_0$ admet un voisinage analytique isomorphe \`a un ouvert analytique de $\hilt\times X^{n-3}$, et dans cette identification le \mor $\pi_{n2}$ s'\'ecrit $\pi_{32}\times id$. On applique le lemme \ref{lema310} pour conclure.

Dans le second cas soit $Z_0=(Z,Z')\in D\subset\hilx\times\dhild$. Le sch\'ema $Z$ est la r\'eunion disjointe $Z'\amalg Z''\amalg x_5\amalg\ldots\amalg x_n$, o\`u $Z''\in\dhild$ et les points de $\supp Z', \supp Z''$ et $x_i$ sont deux \`a deux disjoints. Le point $Z$ admet un voisinage analytique dans $\hilx$ isomorphe \`a un voisinage analytique du point $(Z',Z'',x_5,\cdots,x_n)$ dans $\hild\times\hild\times X^{n-4}$. Dans cette identification, le ferm\'e $D\subset\hilx\times\dhild$ co\"{\i}ncide avec le ferm\'e $D\subset(\hild\times\hild \times X^{n-4})\times\dhild$ des points $((Z',Z'',x_5,\cdots,x_n),Z''')$ tels que $Z'=Z'''$. Alors $D$ est isomorphe localement avec $\dhild\times\hild \times X^{n-4}$, donc $D$ est lisse. \hfill $\Box$

\begin{rem}
\label{remarca63,5}
{\rm On ne peut pas affirmer l'existence d'une application $\tmu:\Sigma\to X$, puisque la vari\'et\'e $\hilxnd$ n'est pas forc\'ement lisse. On ne peut pas appliquer la proposition II 8.23(b) de \cite{Hart} pour d\'emontrer la normalit\'e de la vari\'et\'e $\den$.}
\end{rem}

Les lemmes \ref{lema19} et  \ref{lema20} admettent les g\'e\-n\'e\-ra\-li\-sa\-tions sui\-vantes, dont la d\'e\-mon\-stra\-tion est com\-pl\`ete\-ment a\-na\-logue:

\begin{lema}
\label{lema64}
Le morphisme 
$\elden\to\mu^*\eld\vert_{\dhilxs}$
sur $\hilxs$  de la suite exacte (\ref{ecuatia000}) se prolonge \`a un morphisme surjectif sur $\tihilx$:
$$\elden\to\pi_*\tmu^*\eld.$$ 
\end{lema} 

\begin{lema}
\label{lema65}
Le morphisme 
$\dedoin\to\mu^*\eld\vert_{\dhilxs}$
sur $\hilxs$  de la suite exacte (\ref{ecuatia000}) se prolonge \`a un morphisme surjectif sur $\tihilx$:
\begin{equation*}
\dedoin\to\pi_*\tmu^*\eld.
\end{equation*}
 
\end{lema} 

\vs Le faisceau $\tmu^*\eld$ est \inve sur le grand ouvert $\tisigma=\pi^{-1}(\tihilx)$ (dont le compl\'ementaire est de codimension sup\'erieure ou \'egale \`a $2$) dans la vari\'et\'e normale $\Sigma$. Il se prolonge \`a un faisceau \inve $\widetilde{\eld}$ sur $\Sigma$.
On note $j:\tihilx\to\hilx$ l'inclusion canonique. On est en mesure d'\'enoncer

\begin{teor}
\label{teorema66}
La suite exacte (\ref{ecuatia000}) de faisceaux sur $\hilxs$ se prolonge \`a une suite exacte sur $\tihilx$:
\begin{equation}
\label{ecuatia628}
0\to\esdoi\elen\to j_*\dedoin\oplus \elden\to \pi_*\widetilde{\eld}\to 0.
\end{equation}
\end{teor}

\vs\par {\bf Preuve:}

La preuve reprend \`a l'identique la preuve du \th \ref{teorema13}, pour les inclusions $\hilxs\stackrel{k}{\to}\tihilx\stackrel{j}{\to}\hilx$. On utilise l'annulation 
\begin{equation}
\label{ecuatia630}
\er^1j_*(\esdoi\elen)=0,
\end{equation}
qui vient du fait que $j:\tihilx\to\hilx$ est l'inclusion d'un ouvert dont le compl\'ementaire est de codimension $3$ dans la vari\'et\'e lisse $\hilx$. \hfill $\Box$
%%
%%On note $i:\hilxs\to\tihilx$ l'inclusion canonique. De la m\^eme mani\`ere que dans la d\'emonstration du \th \ref{teorema13} on applique le foncteur  $i_*$ \`a la suite exacte (\ref{ecuatia000}), on utilise le fait que les faisceaux $\esdoi\elen, \dedoin$ et $\elden$ sont localement libres sur $\tihilx$, on utilise les lemmes \ref{lema64} et \ref{lema65} et on obtient une suite exacte sur $\hilx$:
%%$$0\to\esdoi\elen\to j_*\dedoin\oplus\elden\to \pi_*\tmu^*\eld\to 0.$$
%%On applique le foncteur exacte \`a gauche  $j_*$ \`a cette suite exacte. Comme 
%%les faisceaux $\esdoi\elen$ et $\elden$ sont localement libres sur $\hilx$, on obtient une suite exacte sur $\hilx$:
%%\begin{equation}
%%\label{ecuatia629}
%%0\to\esdoi\elen\to j_*\dedoin\oplus\elden\to j_*\pi_*\tmu^*\eld\to\er^1j_*(\esdoi\elen)
%%\end{equation}
%%
%%
%%Consid\'erons le diagramme commutatif
%%\begin{equation}
%%\label{ecuatia631}
%%{\diagram
%%\tisigma\rto^{k}\dto_{\pi}&\Sigma\dto^{\pi}\\
%%\tihilx\rto^{j}&\hilx,
%%\enddiagram}
%%\end{equation}
%%o\`u $k$ est l'inclusion canonique. On a
%%\begin{equation}
%%\label{ecuatia632}
%%j_*\pi_*\tmu^*\eld=\pi_*k_*\tmu^*\eld=\pi_*k_*k^*\tield=\pi_*\tield,
%%\end{equation}
%%en appliquant, dans l'ordre, la commutativit\'e du diagramme (\ref{ecuatia631}), la d\'efinition de $\tield$, le fait que $\tield$ est localement libre et le fait que $k$ est l'inclusion d'un grand ouvert.
%%En rempla\c{c}ant dans la suite exacte (\ref{ecuatia629}) les relations (\ref{ecuatia630}), (\ref{ecuatia632}) on obtient la suite (\ref{ecuatia628}) de l'\'enonc\'e. \hfill $\Box$

\vs La remarque suivante r\'esume les r\'esultats techniques suffisants pour une g\'e\-n\'e\-ra\-li\-sa\-tion des th\'e\-o\-r\`emes \ref{teorema11} (n=2) et \ref{teor24} (n=3) \`a un r\'esultat valable pour $n$ g\'en\'eral. On rappelle que les sch\'emas $\hilxnd$ respectivement $\den$ ont \'et\'es d\'efinis dans le diagramme (\ref{ecuatia0}), respectivement (\ref{ecuatia600}). On rappelle la d\'efinition suivante de \cite{KKMSD}, chap. I,\S 3:

\begin{defi}
\label{def66,5}
Un sch\'ema $X$ est \`a singularit\'es rationnelles s'il est normal et si, pour $f:Z\to X$ une r\'esolution des singularit\'es de $X$, l'une des conditions \'equivalentes  suivantes est v\'erifi\'ee:

a) $\er^if_*\maO_Z=0$ pour $i>0$.

b) $X$ est Cohen-Macaulay et $f_*\omega_Z=\omega_X$

(o\`u $\omega_Z$ est le faisceau dualisant sur $Z$ et $\omega_X$ le faisceau dualisant sur $X$).
\end{defi}

\begin{rem}
\label{prop67}
{\rm Les th\'eor\`emes \ref{teorema11} et \ref{teor24} sugg\`erent pour $n\ge 2$ l'isomorphisme:
$$\has^*(\hilx,\esdoi\elen)=\es^{n-2}\has^*(X,O)\otimes\esdoi\has^*(X,L)\bigoplus(\es^{n-1}\has^*(X,\maO)/{\es^{n-2}\has^*(X,\maO)})\otimes\has^*(X,\eld).$$
Cette affirmation r\'esulte des hypoth\`eses suivantes:
%Soit $n\ge 2$. Si:
\begin{eqnarray}
\label{ecuatia633}
&&\bullet {\rm\ \ le\ \ faisceau\ \ }\dedoin{\rm\ \ est\ \ localement\ \ libre\ \ sur \ \ }\hilx\\
\label{ecuatia634}
&&\bullet {\rm \ \ le \ \ \ \mor\ \ }\pind:\hilxnd\to\hilx{\rm\ \ satisfait\ \ }\er^q\pinds(\pnds\dedoi)=0{\rm\ \ pour\ \ }q>0\\
\label{ecuatia635}
&&\bullet {\rm \ \ le \ \ \mor\ \ }\pind:\den\to\hilx{\rm\ \ satisfait\ \ }\er^q\pinds(\pnds\dedoi)=0{\rm\ \ pour\ \ }q>0\\
\label{ecuatia636}
&&\bullet {\rm\ \ le\ \ sch\acute{e}ma\ \ }\hilxnd{\rm \ \ est\ \ \grave{a}\ \ singularit\acute{e}s \ \ rationnelles\ \ }\\
\label{ecuatia637}
&&\bullet {\rm\ \ le\ \ sch\acute{e}ma\ \ }\den{\rm \ \ est\ \ \grave{a}\ \ singularit\acute{e}s \ \ rationnelles.\ \ }
\end{eqnarray}}
%alors:
%$$\has^*(\hilx,\esdoi\elen)=\es^n\has^*(X,L)\bigoplus(\es^{n-1}\has^*(X,\maO)/{\es^{n-2}\has^*(X,\maO)})\otimes\has^*(X,\eld).$$}
\end{rem}

\vs\par{\bf Preuve de la remarque \ref{prop67}:}

La relation (\ref{ecuatia633}) permet de remplacer $j_*\dedoin$ par $\dedoin$ dans la suite (\ref{ecuatia628}). On commence par d\'emontrer l'analogue de la proposition \ref{prop22}:
\begin{lema}
\label{lema68}
Dans les hypoth\`eses (\ref{ecuatia633}), (\ref{ecuatia634}) et (\ref{ecuatia636}) on a:
$$\has^*(\hilx,\dedoin)=\es^{n-2}\has^*(X,\maO)\otimes\esdoi\has^*(X,L).$$
\end{lema}

\vs\par{\bf Preuve du lemme:}

La relation (\ref{ecuatia634}) implique 
\begin{equation}
\label{ecuatia637,5} 
\has^*(\hilx,\dedoin)=\has^*(\hilxnd,\pnds\dedoi).
\end{equation}
On consid\`ere le \mor $\tp_{n2}=(a,\pnd):\hilxnd\to \es^{n-2}X\times\hild$, o\`u le \mor $a$ associe au couple $(Z,Z')\in\hilx\times\hild$ le point $HC(Z)-HC(Z')$, $HC$ \'etant le \mor de Hilbert-Chow. On consid\`ere ensuite une r\'esolution des singularit\'es $r:\hilxndp\to\hilxnd$. Le \mor $\tp_{n2}$ est birationnel, donc $\tp'_{n2}=\tp_{n2}\circ r$ est 
une r\'esolution des singularit\'es. Par \cite{Bou}, 
$\es^{n-2}X\times\hild$ est \`a singularit\'es rationnelles, donc
\begin{equation}
\label{ecuatia638} 
\er^q\tp'_{n2*}\maO_{\hilxndp}=\begin{cases}0 \mbox{ \ \ \ si \ \ \ } q>0\cr \maO_{\es^{n-2}X\times\hild} \mbox{ \ \ \ si \ \ \ } q=0.\cr\end{cases}
\end{equation}
L'hypoth\`ese (\ref{ecuatia636})  implique
\begin{equation}
\label{ecuatia639} 
\er^qr_{*}\maO_{\hilxndp}=\begin{cases}0 \mbox{ \ \ \ si \ \ \ } q>0\cr \maO_{\es^{n-2}X\times\hild} \mbox{ \ \ \ si \ \ \ } q=0.\cr\end{cases}
\end{equation}

Les relations (\ref{ecuatia638}), (\ref{ecuatia639}) et la suite spectrale de Leray pour $\er\tp'_{n2*}=\er p_{n2*}\circ\er r_*$ impliquent:
$$
\er^qp_{n2*}\maO_{\hilxnd}=\begin{cases}0 \mbox{ \ \ \ si \ \ \ } q>0\cr \maO_{\es^{n-2}X\times\hild} \mbox{ \ \ \ si \ \ \ } q=0\cr\end{cases}$$
d'o\`u:
\begin{equation}
\label{ecuatia640} 
\has^*(\hilxnd,\pnds\dedoi)=\has^*(\es^{n-2}X\times\hild,\maO\boxtimes\dedoi)=\es^{n-2}\has^*(X,\maO)\otimes\esdoi\has^*(X,L).
\end{equation}
Les relations (\ref{ecuatia637,5}) et (\ref{ecuatia640}) impliquent le r\'esultat. \hfill $\Box$

Le lemme suivant g\'en\'eralise la proposition \ref{prop23}:
\begin{lema}
\label{lema69}
Dans les hypoth\`eses (\ref{ecuatia635}) et (\ref{ecuatia637}) on a:
$$\has^*(\hilx,\tield)=\es^{n-2}\has^*(X,\maO)\otimes\has^*(X,\eld).
$$
\end{lema}

\vs\par{\bf Preuve:}

On consid\`ere le \mor $\nu=(a,\tmu):\Sigma\to \es^{n-2}X\times X$, o\`u $a$ est le \mor $\Sigma\to  \es^{n-2}X$ donn\'e par $a(s)=HC(\pi(s))-2\tmu(s).$
En analogie avec la preuve de la proposition \ref{prop23}, il suffit de d\'emontrer:
\begin{lema}
\label{lema70}
On a:
\begin{equation}
\label{ecuatia680}
{\begin{array}{ccc}
   \er^q\nu_*\maO_{\Sigma}&=&\begin{cases}0 \mbox{ \ \ \ si \ \ \ } q>0\cr \maO_{\es^{n-2}X\times X} \mbox{ \ \ \ si \ \ \ } q=0.\cr\end{cases}
 \end{array}}
\end{equation}
\end{lema}

\vs\par{\bf Preuve du lemme \ref{lema70}:}

Soit $r:D'\to\den$ une r\'esolution des singularit\'es. Consid\'erons le diagramme:
$$\diagram
D'\rto^{r}&\den\rto^{\tnu}\dto_{\tpi_{n2}}&\es^{n-2}X\times\dhild\dto^{b}\\
&\Sigma\rto^{\nu}\dto_{\pi}&\es^{n-2}X\times X\\
&\dhilt&
\enddiagram$$
o\`u $\tnu=(a\circ\pi_{n2},p_{n2})$.
En analogie avec la preuve du lemme \ref{lema400}, l'hypoth\`ese (\ref{ecuatia635}) permet de r\'eduire le probl\`eme (\ref{ecuatia680}) au \mor $\den\to\es^{n-2}X\times X$, l'hypoth\`ese (\ref{ecuatia637}) nous permet encore  de nous r\'eduire au \mor $D'\to\es^{n-2}X\times X$, et le fait que $\es^{n-2}X\times\dhild$ est \`a singularit\'es rationnelles 
r\'eduit le probl\`eme au \mor $\es^{n-2}X\times X\to\es^{n-2}X\times X$, fibration en $\proj_1$. \hfill $\Box$

\vs On vient de d\'emontrer les lemmes \ref{lema68} et \ref{lema69}, les analogues des propositions \ref{prop22} et \ref{prop23}. \`A partir de l\`a, la d\'emonstration de la proposition \ref{prop67}  recopie identiquement  la preuve du  \th \ref{teor24}. \hfill $\Box$

\end{document}